\documentclass[dvips]{amsart2000}
\usepackage{epsfig}

\newtheorem{thm}{Theorem}[section]
\newtheorem{cor}[thm]{Corollary}
\newtheorem{lem}[thm]{Lemma}
\newtheorem{prop}[thm]{Proposition}
\theoremstyle{definition}

\newtheorem{rem}[thm]{Remark}
\numberwithin{equation}{section}
\newcommand{\norm}[1]{\left\Vert#1\right\Vert}

\newcommand{\eps}{\varepsilon}

\newcommand{\R}{\mathbb{R}}
\newcommand{\N}{\mathbb{N}}
\newcommand{\Z}{\mathbb{Z}}
\newcommand{\C}{\mathbb{C}}

\newcommand{\dis}{\displaystyle}
\newcommand{\Notequiv}{/\kern-.6em\hbox{$\equiv$} }

\begin{document}

\title[Small polynomials with integer coefficients]
{Small polynomials with integer coefficients}%
\author{Igor E. Pritsker}%

\address{Department of Mathematics, 401 Mathematical Sciences, Oklahoma State
University, Stillwater, OK 74078-1058, U.S.A.}%
\email{igor@math.okstate.edu}

\thanks{Research supported in part by the National Science Foundation
grant DMS-9996410.}%



\begin{abstract}

We study the problem of minimizing the supremum norm, on a segment
of the real line or on a compact set in the plane, by polynomials
with integer coefficients. The extremal polynomials are naturally
called integer Chebyshev polynomials. Their factors, zero
distribution and asymptotics are the main subjects of this paper.
In particular, we show that the integer Chebyshev polynomials for
any infinite subset of the real line must have infinitely many
distinct factors, which answers a question of Borwein and
Erd\'{e}lyi. Furthermore, it is proved that the accumulation set
for their zeros must be of positive capacity in this case.

We also find the first nontrivial examples of explicit integer
Chebyshev constants for certain classes of lemniscates. Since it
is rarely possible to obtain an exact value of integer Chebyshev
constant, good estimates are of special importance.

Introducing the methods of weighted potential theory, we
generalize and improve the Hilbert-Fekete upper bound for integer
Chebyshev constant. These methods also give bounds for the
multiplicities of factors of integer Chebyshev polynomials, and
lower bounds for integer Chebyshev constant. Moreover, all the
mentioned bounds can be found numerically, by using various
extremal point techniques, such as weighted Leja points algorithm.
Applying our results in the classical case of the segment $[0,1]$,
we improve the known bounds for the integer Chebyshev constant and
the multiplicities of factors of the integer Chebyshev
polynomials.

\end{abstract}

\maketitle


\section{Integer Chebyshev problem: History and new results}

Define the uniform (sup) norm on a compact set  $E \subset {\C}$
by
$$\norm{f}_{E} := \sup_{z \in E} |f(z)|.$$
The primary goal of this paper is the study of polynomials with
integer coefficients that minimize the sup norm on the set $E$.
In particular, we consider the asymptotic behavior of these
polynomials and of their zeros. Let ${\mathcal P}_n ( {\C})$
and ${\mathcal P}_n ( {\Z})$ be the classes of algebraic
polynomials of degree at most $n$, respectively with complex
and with integer coefficients. The problem of minimizing the
uniform norm on $E$ by {\it monic} polynomials from ${\mathcal
P}_n ( {\C})$ is well known as the Chebyshev problem (see
\cite{BE95},  \cite{Ri90}, \cite{Ts75}, \cite{Go69}, etc.) In
the classical case $E=[-1,1]$, the explicit solution of this
problem is given by the monic Chebyshev polynomial of degree
$n$:
$$T_n (x) := 2^{1 -n} \cos (n \arccos x), \quad n \in \N.$$
Using a change of variable, we can immediately extend this to an
arbitrary interval $[a,b] \subset {\R}$, so that
$$t_n (x) := \left( \frac{b -a}{2} \right)^n T_n \left(
\frac{2x -a -b}{b -a} \right)$$ is a monic polynomial with real
coefficients and the smallest uniform norm on $[a,b]$ among all
{\it monic} polynomials from ${\mathcal P}_n ( {\C} )$. In fact,
\begin{equation} \label{1.1}
\| t_n \|_{[a,b]} = 2 \left( \frac{b -a}{4} \right)^n, \quad n
\in {\N},
\end{equation}
and we find that the {\it Chebyshev constant} for $[a,b]$ is
given by
\begin{equation} \label{1.2}
t_{\C}([a,b]) := \lim_{n \rightarrow \infty} \| t_n
\|_{[a,b]}^{1/n} = \frac{b -a}{4}.
\end{equation}
The Chebyshev constant of an arbitrary compact set $E \subset
{\C}$ is defined in a similar fashion:
\begin{equation} \label{1.3}
t_{\C}(E) := \lim_{n \rightarrow \infty} \| t_n \|_{E}^{1/n},
\end{equation}
where $t_n$ is the Chebyshev polynomial of degree $n$ on $E$.
It is known that $t_{\C}(E)$ is equal to the transfinite
diameter and the logarithmic capacity $cap(E)$ of the set $E$
(cf. \cite[pp. 71-75]{Ts75}, \cite{Go69} and \cite{Ra95} for the
definitions and background material).

One may notice that the Chebyshev polynomials on the interval
$[-2,2]$ have integer coefficients. The  roots of the $n$-th
Chebyshev polynomial on $[-2,2]$ are
\begin{equation} \label{1.4}
x_k=2\cos\frac{(2k-1)\pi}{2n}, \quad k=1,\ldots,n.
\end{equation}
A remarkable result of Kronecker \cite{Kr1857} states that {\it
any} complete set of conjugate algebraic integers, i.e., roots of
a monic irreducible polynomial over $\Z$, all contained in
$[-2,2]$, must belong to one of the sets (\ref{1.4}) for some
$n\in \N$. Thus we have an exhaustive description of all complete
sets of conjugate algebraic integers in $[-2,2]$, which indicates
that there are infinitely many such sets in this interval. In
fact, Kronecker first proved in \cite{Kr1857} that any complete
set of conjugates on the unit circle $\{|z|=1\}$ must be a subset
of the roots of unity, and then deduced the above result by using
the transformation $x=z+1/z.$ It is difficult to obtain such a
complete characterization when $[-2,2]$ is replaced by a more
general set, but one can extract substantial amount of interesting
information from the study of {\it integer Chebyshev problem}.

An {\it integer Chebyshev polynomial} $Q_n \in {\mathcal P}_n (
{\Z})$ for a compact set $E \subset \C$ is defined by
\begin{equation} \label{1.5}
\| Q_n \|_{E} = \inf_{0 \not\equiv P_n \in{\mathcal P}_n (
{\Z})} \| P_n \|_{E},
\end{equation}
where the $\inf$ is taken over all polynomials from ${\mathcal
P}_n ( {\Z})$, which are not identically zero. Further, the
{\it integer Chebyshev constant} (or integer transfinite
diameter) for $E$ is given by
\begin{equation} \label{1.6}
t_{\Z}(E) := \lim_{n \rightarrow \infty} \| Q_n \|_{E}^{1/n}.
\end{equation}
The existence of the limit in (\ref{1.6}) follows by the same
argument as for (\ref{1.3}), which may be found in \cite{Go69} or
\cite{Ts75}. Note that, for any $P_n \in{\mathcal P}_n ( {\Z})$,
$$ \| P_n \|_{E} = \| P_n \|_{E^*},$$
where $E^*:=E \cup \{z:{\bar z} \in E\}$, because $P_n$ has
real coefficients. Thus the integer Chebyshev problem on a
compact set $E$ is equivalent to that on $E^*$, and we can
assume that $E$ is symmetric with respect to the real axis
($\R$-symmetric) without any loss of generality.

One may readily observe that if $E=[a,b]$ and $b -a \geq 4$,
then $Q_n (x) \equiv 1,\ n \in {\N},$ by (\ref{1.1}) and
(\ref{1.6}), so that
\begin{equation} \label{1.7}
t_{\Z}([a,b]) = 1, \quad b-a \ge 4.
\end{equation}
On the other hand, we obtain directly from the definition and
(\ref{1.2}) that
\begin{equation} \label{1.8}
\frac{b -a}{4} = t_{\C}([a,b]) \leq t_{\Z}([a,b]),\quad b-a<4.
\end{equation}
Hilbert \cite{Hi1894} proved an important upper bound
\begin{equation} \label{1.9}
t_{\Z}([a,b]) \leq \sqrt{\frac{b-a}{4}},
\end{equation}
by using Legendre polynomials and Minkowski theorem on the
integer lattice points in a convex body. Actually, he worked
with $L_2$ norm on $[a,b]$, but this gives the same $n$-th root
behavior as for $L_{\infty}$ norm in (\ref{1.6}).

With the help of Hilbert's result (\ref{1.9}), Schur and Polya
(see \cite{Sch18}) showed that any interval $[a,b] \subset \R$, of
length less than 4, can contain only finitely many complete sets
of conjugate algebraic integers. Thus one may be able to
explicitly find those polynomials with integer coefficients and
all roots in $[a,b],\ b-a<4$. These results were generalized to
the case of an arbitrary compact set $E \subset \C$ by Fekete
\cite{Fe23}, who developed a new analytic setting for the problem,
by introducing the transfinite diameter of $E$ and showing that it
is equal to $t_{\C}(E).$ Both quantities were later proved to be
equal to the logarithmic capacity $cap(E)$, by Szeg\H{o}
\cite{Sz24}. Therefore we state the result of Fekete as follows:
\begin{equation} \label{1.10}
t_{\Z}(E) \leq \sqrt{t_{\C}(E)} = \sqrt{cap(E)},
\end{equation}
where $E$ is $\R$-symmetric. It contains Hilbert's estimate
(\ref{1.9}) as a special case, since $t_{\C}([a,b])=(b-a)/4$ by
(\ref{1.2}). Using the same argument as in \cite{Sch18}, Fekete
concluded by (\ref{1.10}) that there are only finitely many
complete sets of conjugate algebraic integers in any compact set
$E$, satisfying $cap(E)<1$. These ideas found many applications,
but we only discuss here the developments that are closely related
to the subject of this paper. Fekete and Szeg\H{o} \cite{FS55}
showed that any open neighborhood of the set $E$, which is
symmetric in real axis and has $cap(E)=1$, must contain infinitely
many complete sets of conjugates. Robinson \cite{Ro62} proved that
any interval of length greater than 4 carries infinitely many
complete sets of conjugates. But the case of intervals of length
exactly 4, or sets of capacity 1, in general, remains open (for
further references, see \cite{Ro64}, \cite{Ro69}, etc.)

The following useful observation on the asymptotic sharpness for
the estimates (\ref{1.9}) of Hilbert and (\ref{1.10}) of Fekete
is due to Trigub \cite{Tr71}.

\begin{rem} \label{rem1.1}
For the sequence of the intervals $I_m:=[1/(m+4),1/m]$, we have
$$ t_{\Z}(I_m) > \frac{1}{m+2},$$
so that
$$ \lim_{m \to \infty} \left( t_{\Z}(I_m) -
\sqrt{\frac{|I_m|}{4}} \right) = 0.$$
\end{rem}
We include a proof of this fact, due to a relative
inaccessibility of the original paper \cite{Tr71}.

The value $t_{\Z}([a,b])$ is not known for any segment $[a,b],\
b-a<4.$ This represents a difficult open problem, as can be
seen from the study of the classical case $E=[0,1]$, which is
considered below. From a more general point of view, we are
able to find the exact value of $t_{\Z}(E)$ only for a special
class of compact sets, namely for lemniscates. Note that if
$cap(E) \ge 1$ then the problem is trivial, because
$\norm{P_n}_E \ge (cap(E))^n$ for any $P_n \in {\mathcal P}_n (
{\Z})$ of exact degree $n$ (cf. \cite[p. 155]{Ra95}). This
implies that
$$t_{\Z}(E) = 1, \quad \mbox{ if } cap(E) \ge 1.$$

\begin{prop} \label{prop1.2}
Let
\begin{equation} \label{1.11}
V_m(z):=a_mz^m+\ldots+a_0 \in {\mathcal P}_m({\Z}),\ a_m \neq 0.
\end{equation}
Then we have for the lemniscate
\begin{equation} \label{1.12}
L_r:=\{z:|V_m(z)|=r\}, \quad 0 \le r < 1,
\end{equation}
that
\begin{equation} \label{1.13}
(r/|a_m|)^{1/m} \leq t_{\Z}(L_r) \leq r^{1/m}.
\end{equation}
\end{prop}
This gives an immediate corollary.

\begin{cor} \label{cor1.3}
If $V_m(z)$ of (\ref{1.11}) is monic, then
\begin{equation} \label{1.14}
t_{\Z}(L_r) = r^{1/m},
\end{equation}
where $L_r$ is defined in (\ref{1.12}). Furthermore, $(V_m)^k$ is
an integer Chebyshev polynomial of degree $km,\ k \in \N.$
\end{cor}
One may notice that $t_{\Z}(L_r)=t_{\C}(L_r)=cap(L_r)$ (see
\cite[p. 135]{Ra95}) in Corollary \ref{cor1.3}. However, the
following result is more interesting.

\begin{thm} \label{thm1.4}
Suppose that the polynomial $V_m(z)$ of (\ref{1.11}) is
irreducible over integers and that $L_r$ of (\ref{1.12}) satisfies
$0 \le r \le 1/|a_m|$. Then
\begin{equation} \label{1.15}
t_{\Z}(L_r) = r^{1/m},
\end{equation}
and $(V_m)^k$ is an integer Chebyshev polynomial of degree $km,\ k
\in \N.$
\end{thm}
Observe that $t_{\Z}(L_r) \neq t_{\C}(L_r) = cap(L_r) =
(r/|a_m|)^{1/m}$ in this case (cf. \cite[p. 135]{Ra95}).

A deeper insight into the nature of integer Chebyshev constant and
properties of the asymptotically extremal polynomials for integer
Chebyshev problem can be found in the study of this problem for
$E=[0,1]$. It was initiated by Gelfond and Schnirelman, who
discovered an elegant connection with the distribution of prime
numbers (see \cite{GS36} and Gelfond's comments in \cite[pp.
285--288]{Cheb44}). Their argument shows that if $ t_{\Z}([0,1]) =
1/e,$ then the Prime Number Theorem follows. Unfortunately, $
t_{\Z}([0,1]) > 1/e,$ as we shall see below. One can find a nice
exposition of this and related topics in Montgomery \cite[Ch.
10]{Mo94} (also see Chudnovsky \cite{Ch83}). Let ${\mathcal F}_n
\subset {\mathcal P}_n({\Z})$ be the set of irreducible over $\Z$
polynomials, of exact degree $n$, that have all their zeros in
$[0,1].$ Define
\begin{equation} \label{1.16}
s:=\liminf_{n \to \infty \atop F_n \in {\mathcal F}_n}
c_n^{1/n},
\end{equation}
where $F_n=c_nx^n+\ldots.$ Then
\begin{equation} \label{1.17}
t_{\Z}([0,1]) \ge 1/s,
\end{equation}
which is the content of Theorem 2 in \cite[p. 182]{Mo94}. In fact,
Montgomery conjectured that {\it equality} holds in (\ref{1.17}),
but this remains open (essentially the same conjecture was also
made in \cite[p. 90]{Ch83}). One may try to construct various
sequences of polynomials $F_n \in {\mathcal F}_n,\ n \in \N,$ to
obtain lower bounds for $t_{\Z}([0,1])$ from (\ref{1.17}). A few
of such sequences have been devised (cf. \cite{Mo94} and
\cite{Ch83}), with the best known being the Gorshkov sequence of
polynomials. It was originally found by Gorshkov in \cite{Gor59},
and rediscovered by Wirsing \cite{Mo94} and others. These
polynomials arise as the numerators in the sequence of iterates of
the rational function
$$ u(x)=\frac{x(1-x)}{1-3x(1-x)},$$
and they give the following lower bound:
\begin{equation} \label{1.18}
t_{\Z}([0,1]) \ge 1/s_0=0.420726\ldots
\end{equation}
(see \cite[pp. 183-188]{Mo94}).

The upper bounds for $t_{\Z}([0,1])$ can be obtained from the
very definition of integer Chebyshev constant
(\ref{1.5})-(\ref{1.6}). One may even try to find some low
degree integer Chebyshev polynomials and compute their norms,
to find out that this is quite a nontrivial exercise. It was
noticed in many papers that small polynomials from ${\mathcal
P}_n(\Z),\ n \in \N,$ arise as products of powers of
polynomials from ${\mathcal F}_n,\ k<n.$ Aparicio was the first
to prove this in the following strong form (cf. Theorem 3 in
\cite{Ap88}):

If a sequence $Q_n \in {\mathcal P}_n(\Z),\ n \in \N,$ satisfies
\begin{equation} \label{1.19}
\lim_{n \rightarrow \infty} \| Q_n \|_{[0,1]}^{1/n} =
t_{\Z}([0,1]),
\end{equation}
then
\begin{equation} \label{1.20}
Q_n (x) = (x(1 -x))^{[ \alpha_1 n]} (2x -1)^{[ \alpha_2 n]}
(5x^2 -5x +1)^{[\alpha_3 n]} R_n (x), \quad \mbox{as } n \rightarrow \infty,
\end{equation}
where
\begin{equation} \label{1.21}
\alpha_1 \geq 0.1456,\quad \alpha_2 \geq 0.0166 \quad \mbox{ and }
\quad \alpha_3 \geq 0.0037,
\end{equation}
and $R_n \in {\mathcal P}_n ( {\Z})$, $n \in {\N}$.

This gives a good indication of what might be the asymptotic
structure of the integer Chebyshev polynomials on $[0,1]$ and
other sets. Thus Amoroso \cite{Am90} considered intervals with
rational endpoints, and applied a refinement of Hilbert's approach
in \cite{Hi1894} to the polynomials vanishing with high
multiplicities at the endpoints, to improve upon (\ref{1.9}).
Essentially the same ideas were used by Kashin \cite{Ka91} for
dealing with the symmetric intervals $[-a,a]$, in which case one
should consider polynomials with factors $x^k.$

Borwein and Erd\'{e}lyi \cite{BE96} used numerical optimization
techniques to find small polynomials of the form
\begin{equation} \label{1.22}
Q_n (x) = \prod_{i=1}^k Q_{m_i,i}^{[ \alpha_i n]}(x), \quad
0<\alpha_i<1,\ i=1,\dots,k,
\end{equation}
where $Q_{m_i,i} \in {\mathcal P}_{m_i}({\Z})$ and $\sum_{i=1}^k
\alpha_i m_i =1$. They improved the upper bound for
$t_{\Z}([0,1])$, which triggered a number of numerical studies
on the integer Chebyshev polynomials for $[0,1]$ and other
intervals. Borwein and Erd\'{e}lyi also improved the result of
Aparicio (\ref{1.19})-(\ref{1.21}):
$$ \alpha_1 \ge 0.26,$$
and used this to show that the {\it strict} inequality holds in
(\ref{1.18}). Hence the Gorshkov polynomials do not give the
exact value of $t_{\Z}([0,1])$.

The ideas of Borwein and Erd\'{e}lyi have been developed in the
papers by Flammang \cite{Fl2}, by Flammang, Rhin and Smyth
\cite{FRS97}, and by Habsieger and Salvy \cite{HS97}, to obtain
further numerical improvements in the upper bounds for $t_{\Z}$
on $[0,1]$ and on Farey intervals. In particular, Habsieger and
Salvy computed 75 first integer Chebyshev polynomials for
$[0,1]$ and found the best known upper bound
\begin{equation} \label{1.23}
t_{\Z}([0,1]) \leq 0.42347945.
\end{equation}
Flammang, Rhin and Smyth \cite{FRS97} generalized the approach
of \cite{BE96} to improve the lower bounds in (\ref{1.21})
$$ \alpha_1 \geq 0.264151,\quad \alpha_2 \geq
0.021963 \quad \mbox{ and } \quad \alpha_3 \geq 0.005285, $$ as
well as bounds for six additional factors of the integer
Chebyshev polynomials on $[0,1]$. They also extended the
Gorshkov polynomials technique to the Farey intervals
$[p/q,r/s],$ with $qr-ps=1$, and obtained an interesting
generalization of (\ref{1.18}).

From the above discussion, it is natural to expect that the
integer Chebyshev polynomials for $[0,1]$ are built out of the
factors as in (\ref{1.22}), which is suggested in Montgomery
\cite[p. 182]{Mo94}. In addition, Montgomery proposed to study the
zero distribution of these polynomials, associated measures and
extremal potentials. Potential theory indeed provides powerful
methods for dealing with various extremal problems for
polynomials, which proved to be very effective for classical
Chebyshev polynomials, orthogonal polynomials, etc. It is clear
that the study of zeros for integer Chebyshev polynomials is
essentially equivalent to the study of their factors and
asymptotic behavior. We should note that not all of the zeros of
the integer Chebyshev polynomials for $[0,1]$ actually lie on
$[0,1]$. This was discovered by Habsieger and Salvy \cite{HS97},
who found a factor of an integer Chebyshev polynomial of degree
70, with two pairs of complex conjugate roots.

One might hope that the sequence of the integer Chebyshev
polynomials for $[0,1]$ is composed from products of powers of
a {\it finite } number of irreducible polynomials over $\Z.$
Unfortunately, this is not true as we show by the following
result, answering a question of Borwein and Erd\'{e}lyi (see
\cite{BE96}, Q7).

\begin{thm} \label{thm1.5}
Let $E \subset \R$ be a compact set, $cap(E)<1$, consisting of
infinitely many points. The integer Chebyshev polynomials $Q_{n}$
for $E,\ n \in \N,$ have infinitely many distinct factors with
integer coefficients, as $n \to \infty$.
\end{thm}
It is obvious from the known results that integer Chebyshev
polynomials are completely different from the classical companions
in their ``discrete" nature. However, their zeros cannot be so
isolated, as it might appear.
\begin{thm} \label{thm1.6}
Let $Z$ be the set of accumulation points for the zeros of the
integer Chebyshev polynomials for a compact set $E \subset \R,\
0<cap(E)<1$. Then
\begin{equation} \label{1.24}
cap(Z) > 0.
\end{equation}
\end{thm}
This immediately implies that $Z$ cannot be too small, e.g., it
cannot be a countable set. One might conjecture that the zeros of
the integer Chebyshev polynomials on $[0,1]$ are dense in a
Cantor-type set of positive capacity.

Since the nature of the unknown factors of the integer Chebyshev
polynomials for $[0,1]$ is rather obscure, we may view the
integer Chebyshev polynomials as being of the form
\begin{equation} \label{1.26}
Q_n (x) = \left( \prod_{i=1}^k Q_{m_i,i}^{l_i(n)}(x) \right)
R_n(x), \quad n \in \N,
\end{equation}
where $l_i(n) \in \N,$ $Q_{m_i,i}(x)$ is the known irreducible
factor of degree $m_i,\ i=1,\ldots,k,$ and $R_n(x)$ is the
remainder. Assuming that the limits
\begin{equation} \label{1.27}
\lim_{n \to \infty} \frac{l_i(n)}{n}=:\alpha_i>0, \quad
i=1,\ldots,k,
\end{equation}
exist, at least along a subsequence, we observe that the $n$-th
root of the absolute value of the product in (\ref{1.26})
converges to a fixed ``weight" function, as $n \to \infty$,
locally uniformly in $\C$:
$$ \lim_{n \to \infty} \left( \prod_{i=1}^k |Q_{m_i,i}(x)|^{l_i(n)}
\right)^{1/n} = \prod_{i=1}^k |Q_{m_i,i}(x)|^{\alpha_i},$$ where
$\sum_{i=1}^k \alpha_im_i \le 1.$ Hence, for the purposes of
studying the asymptotic behavior, as $n \to \infty,$ we may regard
$Q_n(x)$ of (\ref{1.26}) as a ``weighted polynomial" and use the
methods of weighted potential theory \cite{ST97}. Following this
idea, we generalize the Hilbert-Fekete upper bound for $t_{\Z}$
and find new lower bounds. We also prove various results on the
multiplicities of factors and zeros of integer Chebyshev
polynomials in the next section. Then we apply the general theory
to the integer Chebyshev problem on $[0,1]$ and obtain substantial
improvements over the previously known results in Section 3.
Section 4 contains a brief outline of the basic facts of weighted
potential theory, used in this paper. All proofs are given in
Section 5.

It must be mentioned that the history of the problem as sketched
here is far from being complete. Integer Chebyshev problem is
closely connected to approximation by polynomials with integer
coefficients (see Ferguson \cite{Fer80} and Trigub \cite{Tr71} for
surveys), which has interesting history of its own. Further
related topics are entire functions with integer coefficients (or
integer valued) (cf. P\'{o}lya \cite{Po15}, \cite{Po22} and
\cite{Po28}, Pisot \cite{Pi42}, \cite{Pi46-1} and \cite{Pi46-2},
and Robinson \cite{Ro68}, \cite{Ro71}, etc.), integer moment
problem (see Barnsley, Bessis and Moussa \cite{BBM79}),
Schur-Siegel trace problem (cf. Schur \cite{Sch18}, Siegel
\cite{Si45}, Smyth \cite{Sm84}, Borwein and Erd\'{e}lyi
\cite{BE96}, etc.) and many others.

\section{Upper and lower bounds for integer Chebyshev constant}

Motivated by the known results on the asymptotic structure of the
integer Chebyshev polynomials, we study the weighted polynomials
$w^n(z)P_n(z)$, where $w(z)$ is a continuous nonnegative function
on a compact $\R$-symmetric set $E \subset \C$ and $P_n \in
{\mathcal P}_n(\Z)$. By analogy with (\ref{1.5})-(\ref{1.6}),
consider the weighted integer Chebyshev polynomials $q_n \in
{\mathcal P}_n(\Z),\ n\in\N,$ such that
$$v_n(E,w):=\norm{w^n q_n}_E = \inf_{0
\not\equiv P_n \in{\mathcal P}_n ( {\Z})} \| w^n P_n \|_{E},$$ and
define the weighted integer Chebyshev constant by
\begin{equation} \label{2.1}
t_{\Z}(E,w) := \lim_{n \rightarrow \infty} \left( v_n(E,w)
\right)^{1/n}.
\end{equation}
The limit in (\ref{2.1}) exists by the following standard
argument. Note that
$$ v_{k+m}(E,w) \le \norm{w^{k+m} q_k q_m}_E \le \norm{w^{k} q_k}_E
\norm{w^{m} q_m}_E = v_{k}(E,w) v_{m}(E,w).$$ If we set $a_n =
\log{v_n(E,w)}$, then
$$a_{k+m} \le a_k + a_m,\quad k,m \in\N.$$
Hence
$$\lim_{n\to\infty} \frac{a_n}{n} = \lim_{n\to\infty} \log \left(
v_n(E,w) \right)^{1/n}$$ exists by Lemma on page 73 of
\cite{Ts75}.

Our first goal is to give an upper bound for $t_{\Z}(E,w)$. It
is possible to generalize the Hilbert-Fekete method for this
purpose, but we also need the concept of the {\it weighted
capacity of $E$}, denoted by $cap(E,w)$ (see \cite{ST97} and a
brief overview of the weighted potential theory in Section 4).

\begin{thm} \label{thm2.1}
Let $E \subset \R$ be a compact set and let $w : E \to
[0,+\infty)$ be a continuous function. Then
\begin{equation} \label{2.2}
t_{\Z}(E,w) \le \sqrt{cap(E,w)}.
\end{equation}
\end{thm}

\begin{rem} \label{rem2.2}
If $w(z)\equiv 1$ on $E$ then $cap(E,1)=cap(E)$, i.e., (\ref{2.2})
reduces to the result of Fekete (\ref{1.10}).
\end{rem}

It is clear from Section 1 that our main applications are
related to the weights of the following type:
\begin{equation} \label{2.3}
w(z) = \left( \prod_{i=1}^k |Q_{m_i,i}(z)|^{\alpha_i}
\right)^{1/(1-\alpha)},
\end{equation}
where factors $Q_{m_i,i} \in {\mathcal P}_{m_i}({\Z})$ have the
form
\begin{equation} \label{2.4}
Q_{m_i,i}(z)=a_i\prod_{j=1}^{m_i}(z-z_{j,i}), \quad a_i\neq0,\
i=1,\dots,k,
\end{equation}
and
\begin{equation} \label{2.5}
\alpha:=\sum_{i=1}^k \alpha_im_i < 1,
\end{equation}
with $0<\alpha_i<1,\ i=1,\dots,k.$ Thus we immediately obtain
an upper bound for the classical (not weighted) integer
Chebyshev constant.

\begin{thm} \label{thm2.3}
Suppose that $E \subset \R$ is a compact set, and that the weight
$w(z)$ satisfies (\ref{2.3})-(\ref{2.5}). Then
\begin{equation} \label{2.6}
t_{\Z}(E) \le \left(cap(E,w)\right)^{(1-\alpha)/2}.
\end{equation}
\end{thm}
Theorem \ref{thm2.3} suggests that we may be able to improve the
results of Hilbert (\ref{1.9}) and of Fekete (\ref{1.10}), by
using (\ref{2.6}) with a proper choice of factors $Q_{m_i,i},\
i=1,\dots,k,$ for the weight $w.$ It is natural to utilize the
known factors of integer Chebyshev polynomials for that purpose.
We shall carry out this program in the next section, and obtain
an improvement of the upper bound (\ref{1.23}).

It is clear that we need an effective method of finding weighted
capacity, in order to make the estimate (\ref{2.6}) practical. For
the ``polynomial-type" weights we are considering here, one can
express $cap(E,w)$ through the regular logarithmic capacity and
Green functions.

\begin{thm} \label{thm2.4}
Let $E \subset \R$ be a compact set, $cap(E)>0$, and let $w(z)$ be
as in (\ref{2.3})-(\ref{2.5}). Then there exists a compact set
$S_w \subset E \setminus \cup_{i=1}^k \{z_{j,i}\}_{j=1}^{m_i},$
such that (\ref{2.6}) holds with
\begin{equation} \label{2.7}
cap(E,w) = \exp\left(\int\log w\, d\mu_w - F_w \right),
\end{equation}
where
\begin{equation} \label{2.8}
F_w = \frac{1}{\alpha-1}\left(\log{cap(S_w)} + \sum_{i=1}^k
\alpha_i\log|a_i| + \sum_{i=1}^k \sum_{j=1}^{m_i} \alpha_i
g_{\Omega}(z_{j,i},\infty) \right)
\end{equation}
and
\begin{equation} \label{2.9}
\mu_w = \frac{1}{1-\alpha}\left(\omega(\infty,\cdot,\Omega) -
\sum_{i=1}^k \sum_{j=1}^{m_i} \alpha_i
\omega(z_{j,i},\cdot,\Omega) \right)
\end{equation}
is the unit positive measure supported on $S_w$. Alternatively,
\begin{equation} \label{2.10}
cap(E,w) = cap(S_w)\exp\left(\int\log w\, d(
\omega(\infty,\cdot,\Omega) + \mu_w) \right).
\end{equation}
Here, $\Omega:=\overline{\C} \setminus S_w,\ g_{\Omega}(z,\xi)$
is the Green function of $\Omega$ with pole at $\xi \in \Omega$,
and $\omega(\xi,\cdot,\Omega)$ is the harmonic measure at $\xi
\in \Omega$ with respect to $\Omega$.
\end{thm}
Note that $\mu_w$ arises as the equilibrium measure in the
weighted energy problem associated with the weight $w$ of
(\ref{2.3})-(\ref{2.5}), and $F_w$ is the {\it modified Robin
constant} for that energy problem (cf. \cite{ST97} and Section
4 of this paper for the details). The measure
$\omega(\infty,\cdot,\Omega)$ is the classical equilibrium
distribution on $S_w$, in the sense of logarithmic potential
theory (see \cite{Ts75}, \cite{Ra95}, etc.)

Using certain information on the asymptotic behavior of integer
Chebyshev polynomials, we can find lower bounds for integer
Chebyshev constant, as below.

\begin{thm} \label{thm2.5}
Suppose that the integer Chebyshev polynomials of a compact set
$E \subset \C,\ cap(E)>0,$ satisfy, along a subsequence of $n
\to \infty,$
\begin{equation} \label{2.11}
Q_n (z) = \left( \prod_{i=1}^k Q_{m_i,i}^{l_i(n)}(z) \right)
R_n(z),\quad \deg Q_n=n,
\end{equation}
where $Q_{m_i,i}(z) \in {\mathcal P}_{m_i}({\Z})$, $l_i(n) \in
\N,$ and the limits
\begin{equation} \label{2.12}
\lim_{n \to \infty} \frac{l_i(n)}{n}=:\alpha_i > 0, \quad
i=1,\ldots,k,
\end{equation}
exist. Then
\begin{equation} \label{2.13}
t_{\Z}(E) \ge e^{(\alpha-1)F_w},
\end{equation}
where $F_w$  is the modified Robin constant for the weight $w$
of (\ref{2.3}) and $\alpha$ is given by (\ref{2.5}).

Moreover, if $E \subset \R$ then
\begin{equation} \label{2.14}
t_{\Z}(E) \ge cap(S_w) \prod_{i=1}^k |a_i|^{\alpha_i} \exp\left(
\sum_{i=1}^k \sum_{j=1}^{m_i} \alpha_i g_{\Omega}(z_{j,i},\infty)
\right),
\end{equation}
in the notations of Theorem \ref{thm2.4}.
\end{thm}
Theorem \ref{thm2.5} is an easy consequence of the results in
weighted potential theory and a simple fact that the leading
coefficient of $R_n(z)$ is at least 1 in absolute value, being a
nonzero integer. It turns out that we can obtain better lower
bounds for $t_{\Z}(E)$, by using rational points. One of the
possible results in this direction is given below. Recall that the
logarithmic potential of a Borel measure $\mu$ is defined by
$$U^{\mu}(z) := \int \log \frac{1}{|z -t|} d \mu(t).$$

\begin{thm} \label{thm2.6}
Assume that the integer Chebyshev polynomials of $E,\ cap(E)>0,$
satisfy (\ref{2.11}) and (\ref{2.12}), where $R_n(\zeta)\neq0$ for
a point $\zeta\in\C$, along a subsequence of $n \to \infty.$ If
$\zeta=(p_1+ip_2)/q$ is a complex rational number in reduced form,
i.e., $gcd(p_1,p_2,q)=1,$ then
\begin{equation} \label{2.15}
t_{\Z}(E) \ge
q^{\alpha-1}\exp\left((\alpha-1)(F_w-U^{\mu_w}(\zeta))\right),
\end{equation}
where $F_w$ is the modified Robin constant and $\mu_w$ is the
weighted equilibrium distribution associated with the weight
$w$ of (\ref{2.3})-(\ref{2.5}). For $\zeta=0$, we set $q=1$ in
(\ref{2.15}).
\end{thm}
Estimate (\ref{2.15}) has interesting applications in the
``opposite" direction, as it can be used to improve the bounds for
the multiplicities of the known factors of integer Chebyshev
polynomials. Thus we can deduce the ``asymptotic structure" result
from the upper bound for $t_{\Z}(E),$ as an immediate corollary of
Theorem \ref{thm2.6}.

\begin{cor} \label{cor2.7}
Suppose that the assumptions of Theorem \ref{thm2.6} are
satisfied  and that
$$t_{\Z}(E) \le M. $$
Then the set of multiplicities $\{\alpha_i\}_{i=1}^k$ must
satisfy
\begin{equation} \label{2.16}
q^{\alpha-1}\exp\left((\alpha-1)(F_w-U^{\mu_w}(\zeta))\right)\le
M.
\end{equation}
\end{cor}
This inequality defines a domain for the possible values of
$\alpha_i,\ i=1,\ldots,k,$ which allows to significantly
improve the known bounds for $\alpha_i$'s.

Another immediate, but nontrivial, consequence of the weighted
potential theory is the following fact.

\begin{prop} \label{prop2.8}
Let $E \subset \C$ be a compact set, $cap(E)>0.$ Suppose that the
integer Chebyshev polynomials for $E$ satisfy (\ref{2.11}) and
(\ref{2.12}), along a subsequence of $n \to \infty.$ Then there
exists $\eps>0$, so that
\begin{equation} \label{2.17}
t_{\Z}(E)=t_{\Z}(E\cup H_{\eps}),
\end{equation}
where $$H_{\eps}=\bigcup_{i=1}^k\bigcup_{j=1}^{m_i}
\{z:|z-z_{j,i}|\le\eps\}.$$
\end{prop}

Perhaps, the most interesting application of our general results,
developed in this section, is the classical case $E=[0,1]$.
Therefore, we concentrate on its study below, to demonstrate the
strength of the method.

\section{Integer Chebyshev problem on $[0,1]$}

We remind that the best known bounds for $t_{\Z}([0,1])$, as
mentioned in (\ref{1.18}) and (\ref{1.23}), are as follows:
$$0.42072638... < t_{\Z}([0,1]) \leq 0.42347945.$$
The above lower bound was believed to be the precise value of
$t_{\Z}([0,1])$, but Borwein and Erd\'{e}lyi \cite{BE96} showed
that there must be the {\it strict} inequality. However, they
did not give a numerical value for the improvement in the lower
bound. Using the general methods of Section 2, based on
weighted potential theory, we show here that

\begin{thm} \label{thm3.1}
$$ 0.4213 < t_{\Z}([0,1]) < 0.4232.$$
\end{thm}

It is convenient for technical reasons to use the symmetry of
$[0,1]$ and the standard change of variable $x(1 -x) \rightarrow
z$, which reduces the integer Chebyshev problem on $[0,1]$ to
that on $[0,1/4]$:
\begin{equation} \label{3.1}
(t_{\Z}([0,1]))^2 = t_{\Z}([0,1/4]).
\end{equation}
Furthermore, we have by Lemmas 1-2 of \cite{HS97} that the
integer Chebyshev polynomials for $[0,1]$ and $[0,1/4]$ are
related by
\begin{equation} \label{3.2}
Q_{2k} (x) = q_k (x(1 -x))
\end{equation}
and
\begin{equation} \label{3.3}
Q_{2k +1} (x) = (1 -2x)q_k (x(1 -x)).
\end{equation}
Hence we can study the integer Chebyshev problem on $[0,1/4]$,
and then return to $[0,1]$ without any loss of information.

Habsieger and Salvy \cite{HS97} give the following list of
known factors of the integer Chebyshev polynomials for $[0,1]$:
\begin{eqnarray*}
A_1(x)=x(1-x),\quad A_2(x)=2x-1,\quad A_3(x)=5x^2-5x+1,\\
A_4(x)=6x^2-6x+1,\quad A_5(x)=29x^4-58x^3+40x^2-11x+1,\\
A_6(x)=(13x^3-20x^2+9x-1)(13x^3-19x^2+8x-1),\\
A_7(x)=(31x^4-63x^3+44x^2-12x+1)(31x^4-61x^3+41x^2-11x+1),\\
A_8(x)=4921x^{10}-24605x^9+53804x^8-67586x^7+53866x^6-28388x^5\\
+9995x^4-2317x^3+338x^2-28x+1.
\end{eqnarray*}
Incidentally, $A_8(x)$ is the ``surprise factor" with four
non-real zeros. Changing the variable to $z=x(1-x)$, we obtain
the following factors for $[0,1/4]$:
\begin{eqnarray} \label{3.4}
Q_{1,1}(z)=z,\quad Q_{1,2}(z)=A^2_2(x)=4z-1,\quad Q_{1,3}(z)=5z-1, \\
Q_{1,4}(z)=6z-1,\quad Q_{2,5}(z)=29z^2-11z+1,\nonumber\\
Q_{3,6}(z)=169z^3-94z^2+17z-1,\nonumber\\
Q_{4,7}(z)=961z^4-712z^3+194z^2-23z+1,\nonumber\\
Q_{5,8}(z)=4921z^5-4594z^4+1697z^3-310z^2+28z-1.\nonumber
\end{eqnarray}
Exactly these factors will be used in the definition of the
weight $w$ of (\ref{2.3}) for the applications of the results
from Section 2. We start with the case of two factors $Q_{1,1}$
and $Q_{1,2}$, vanishing at the endpoints of $[0,1/4]$, where
all the parameters of the corresponding weighted potential
theory can be found explicitly.

\subsection{Two factors on $[0,1/4]$}

Note that if the relative multiplicities for the factors
$A_1(x)=x(1-x)$ and $A_2(x)=2x-1$, in the integer Chebyshev
polynomial on $[0,1]$, are $\alpha_1$ and $\alpha_2$, then the
relative multiplicities for the corresponding factors
$Q_{1,1}(z)=z$ and $Q_{1,2}(z)=4z-1$ on $[0,1/4]$ are
$2\alpha_1$ and $\alpha_2$ (see (\ref{3.2})-(\ref{3.4})). Thus
we define the weight $w$ according to (\ref{2.3}):
\begin{equation} \label{3.5}
w(x)=\left(x^{2\alpha_1}(1-4x)^{\alpha_2}\right)^{1/(1-2\alpha_1-\alpha_2)},
\quad x\in [0,1/4],
\end{equation}
where $\alpha_1,\alpha_2>0$ and $2\alpha_1+\alpha_2<1.$ The
needed quantities of the weighted potential theory are
contained in the following lemma.

\begin{lem} \label{lem3.2}
For the weight $w$ of (\ref{3.5}), we have that $S_w=[a,b] \subset
(0, 1/4)$, with
\begin{equation} \label{3.6}
a := (4\alpha_1^2 - \alpha_2^2 - \sqrt{\Delta} +1)/8 \quad \mbox{
and } \quad b := (4\alpha_1^2 - \alpha_2^2 + \sqrt{\Delta}+1)/8,
\end{equation}
where $\Delta := (1 -(2\alpha_1 + \alpha_2)^2)(1 - (2\alpha_1 -
\alpha_2)^2)$. Furthermore,
\begin{eqnarray} \label{3.7}
\quad F_w&=&\frac{1-\alpha_2}{1-2\alpha_1-\alpha_2}\log{4} -
\log(b-a) -
\frac{4\alpha_1}{1-2\alpha_1-\alpha_2}\log(\sqrt{a}+\sqrt{b})
\\ &-& \frac{2\alpha_2}{1-2\alpha_1-\alpha_2}
\log(\sqrt{1/4-a}+\sqrt{1/4-b}), \nonumber
\end{eqnarray}
\begin{equation} \label{3.8}
d\mu_w(x)=\frac{\sqrt{(x-a)(b-x)}\ dx
}{\pi(1-2\alpha_1-\alpha_2)x(1/4-x)}, \quad x \in [a,b],
\end{equation}
and
\begin{eqnarray} \label{3.9}
F_w-U^{\mu_w}(z) &=& (g_{\Omega} (z, \infty ) - 2\alpha_1( \log
|z| + g_{\Omega} (z,0)) \\ &-& \alpha_2 ( \log |4z -1| +
g_{\Omega} (z, 1/4))) /(1 - 2\alpha_1 - \alpha_2), \nonumber
\end{eqnarray}
where we use the notation of Theorem \ref{thm2.4}.
\end{lem}
Note that the function in (\ref{3.9}) is continuous in ${\C}$ and
harmonic in ${\C} \setminus [a,b]$. The weighted capacity
$cap([0,1/4],w)$ is found from (\ref{2.7}) or (\ref{2.10}), with
the help of (\ref{3.6})-(\ref{3.8}). Obviously, $cap([0,1/4],w)$
is a function of two variables $\alpha_1$ and $\alpha_2$, defined
on the triangle $T:=\{\alpha_1,\alpha_2>0:
2\alpha_1+\alpha_2<1\}.$ Thus we obtain from Theorem \ref{thm2.3},
via computation, that
$$t_{\Z}([0,1/4])\le\inf_{\alpha_1,\alpha_2\in T}
\left(cap([0,1/4],w)\right)^{(1-2\alpha_1-\alpha_2)/2} \approx
0.18043338.$$ The bound is attained for $\alpha_1 \approx
0.290447$ and $\alpha_2\approx 0.09$, which matches the result
of \cite[p. 906]{Am90}. But this upper bound is greater than the
one in (\ref{1.23}), by (\ref{3.1}), so that it is not
interesting for us.

We can also apply Theorem \ref{thm2.6} here, with $\zeta_1=0$ and
$\zeta_2=1/4$, because these zeros are absorbed by the weight $w$.
Hence we have two simultaneous lower bounds
$$t_{\Z}([0,1/4])>l_1(\alpha_1,\alpha_2):=\exp\left(
(2\alpha_1+\alpha_2-1)(F_w-U^{\mu_w}(0))\right)$$ and
$$t_{\Z}([0,1/4])>l_2(\alpha_1,\alpha_2):=4^{2\alpha_1+\alpha_2-1}
\exp\left( (2\alpha_1+\alpha_2-1)(F_w-U^{\mu_w}(1/4))\right),$$
for $\alpha_1,\alpha_2\in T.$ It follows that
$$t_{\Z}([0,1/4])\ge\inf_{\alpha_1,\alpha_2\in T} \max\left(
l_1(\alpha_1,\alpha_2),l_2(\alpha_1,\alpha_2)\right) \approx
0.176056,$$ where the numerical value, attained for $\alpha_1
\approx 0.330333$ and $\alpha_2\approx 0.128$, is found by using
(\ref{3.9}) and computations. Again, this lower bound is weaker
than (\ref{1.18}).

However, the application of Corollary \ref{cor2.7}, with the upper
bound $M$ obtained from (\ref{1.23}) and (\ref{3.1}), gives an
interesting new result (also see \cite{Pr99}). We translate it to
$[0,1]$ setting here.

\begin{thm} \label{thm3.3}
The integer Chebyshev polynomials $\{ Q_n \}_{n =1}^{\infty}$
on $[0,1]$ satisfy
\begin{equation} \label{3.10}
Q_n (x)=(x(1 -x))^{[ \alpha_1 n]} (2x -1)^{[ \alpha_2 n]} R_n
(x), \quad \mbox{as } n \rightarrow \infty,
\end{equation}
where
\begin{equation} \label{3.11}
0.2961 \leq \alpha_1 \leq 0.3634 \quad \mbox{ and } \quad
0.0952 \leq \alpha_2 \leq 0.1767,
\end{equation}
and $R_n \in {\mathcal P}_n ({\Z})$, $n \in {\N}$. Furthermore,
the pair $(\alpha_1, \alpha_2)$ must belong to the region $G$
pictured below in Figure 1, which is determined by the
inequalities
$$\exp\left((2\alpha_1+\alpha_2-1)(F_w-U^{\mu_w}(0))\right)
< 0.179335$$ and
$$4^{2\alpha_1+\alpha_2-1} \exp\left( (2\alpha_1+\alpha_2-1)
(F_w-U^{\mu_w}(1/4))\right) < 0.179335.$$
\begin{figure}[htb]
\centerline{\psfig{figure=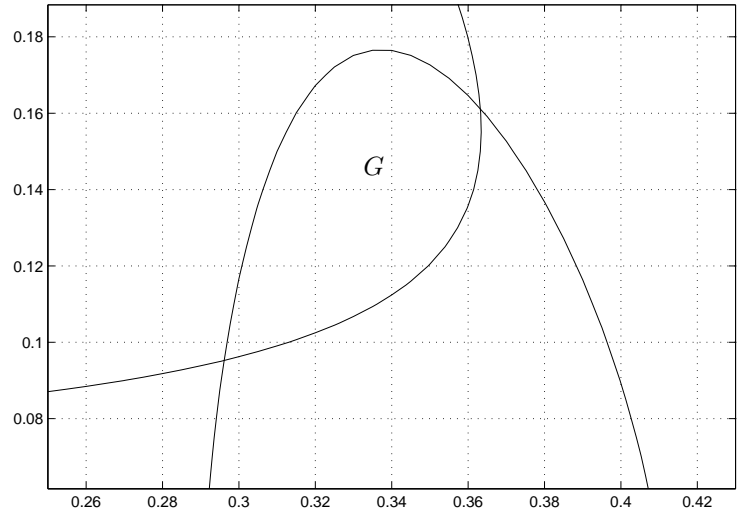,height=3in}} \vskip -2in $G$
\vskip 2in \caption{Region $G$ for $\alpha_1$ and $\alpha_2$.}
\label{fig1}
\end{figure}
\end{thm}
Note that, in addition to improving the previous lower bounds
obtained in \cite{Ap88}, \cite{BE96} and \cite{FRS97},
(\ref{3.11}) also gives the upper bounds for $\alpha_1$ and
$\alpha_2$.

\subsection{Three and more factors on $[0,1/4]$: Numerical approach}

It is natural to expect improvements in the bounds for
$t_{\Z}([0,1/4])$ and for $\alpha_i$'s, if we use three or more
known factors from (\ref{3.4}). There is, however, a substantial
difficulty arising on our way. Although Theorem \ref{thm2.4} can
still produce the needed quantities of weighted potential theory,
it assumes the knowledge of the set $S_w$. In fact, when $w$ is
defined by (\ref{2.3})-(\ref{2.5}), with the help of the factors
\begin{equation} \label{3.12}
Q_{1,1}(z)=z, \quad Q_{1,2}(z)=4z-1 \quad\mbox{and}\quad
Q_{1,3}(z)=5z-1,
\end{equation}
we have that $S_w=[a_1,b_1]\cup [a_2,b_2]$, where
$[a_1,b_1]\subset(0,1/5)$ and $[a_2,b_2]\subset(1/5,1/4)$. But the
endpoints of the intervals $[a_1,b_1]$ and $[a_2,b_2]$ are some
unknown functions of the multiplicities $(\alpha_1,\alpha_2,
\alpha_3)$. The problem becomes even more complicated, if we
consider further factors listed in (\ref{3.4}). Fortunately, we
can apply the numerical methods for finding $S_w$ and the weighted
equilibrium measure $\mu_w$, based on weighted Leja points (see
Section V.1 of \cite{ST97}). Weighted Leja points are easy to
generate numerically, as they are defined by the following simple
recursive procedure. For a general compact set $E$ and $w$ of
(\ref{2.3})-(\ref{2.5}), let $a_0\in E$ be a point such that
\begin{equation} \label{3.13}
|a_0|w(a_0) = \norm{zw(z)}_E.
\end{equation}
Given the points $\{a_i\}_{i=0}^{n-1},$ we define the weighted
Leja polynomial
\begin{equation} \label{3.14}
L_n(z)=\prod_{i=0}^{n-1} (z-a_i),
\end{equation}
so that $a_n$ is found as a point satisfying
\begin{equation} \label{3.15}
w^n(a_n) |L_n(a_n)| = \norm{w^nL_n}_E,\quad a_n\in E.
\end{equation}
Of course, the choice of $a_n$ might not be unique. The
fundamental property of weighted Leja points is that they give a
discrete approximation to the weighted equilibrium measure
$\mu_w,$ corresponding to the weight $w$ on $E,\ cap(E)>0$. This
is stated in the best way by using the weak* convergence of
measures:
\begin{equation} \label{3.16}
\tau_n:=\frac{1}{n+1}\sum_{i=0}^n \delta_{a_i}
\stackrel{*}{\rightarrow} \mu_w,\quad \mbox{as } n\to\infty,
\end{equation}
where $\delta_{a_i}$ is the unit point mass at $a_i,\
i=0,1,2,\ldots$ (see Theorem V.1.1 in \cite{ST97}). Furthermore,
Theorem V.1.2 of \cite{ST97} gives that
\begin{equation} \label{3.17}
\lim_{n\to\infty} w(a_n) |L_n(a_n)|^{1/n} = e^{-F_w}.
\end{equation}
It follows from Theorem III.2.1 and Remark III.2.2 of \cite{ST97}
(also cf. (\ref{4.4})) that $\{a_i\}_{i=0}^{\infty} \subset
S_w^*$, where $S_w^* \supset S_w$ is a compact set defined by
\begin{equation} \label{3.18}
S_w^*:=\{z\in E: U^{\mu_w}(z)-\log w(z) \le F_w \}.
\end{equation}
One can immediately see from (\ref{3.18}) and the form of $w$ in
(\ref{2.3})-(\ref{2.5}), that
\begin{equation} \label{3.19}
S_w^* \subset E\setminus \bigcup_{i=1}^k \{z_{j,i}\}_{j=1}^{m_i}.
\end{equation}
Thus we obtain from (\ref{3.16}), (\ref{3.17}) and the definition
of weak* convergence that
\begin{equation} \label{3.20}
\exp(F_w- U^{\mu_w}(\zeta))=\lim_{n\to\infty} \frac{1}{w(a_n)}
\left|\frac{L_n(\zeta)}{L_n(a_n)}\right|^{1/n},\quad \zeta\in
\C\setminus S_w^*.
\end{equation}
Similarly, we have from (\ref{2.7}), (\ref{3.16}), (\ref{3.17})
and (\ref{3.19}) that
\begin{equation} \label{3.21}
cap(E,w)=\lim_{n\to\infty} w(a_n)
\left(|L_n(a_n)|\prod_{i=0}^{n-1} w(a_i)\right)^{1/n}.
\end{equation}
Equations (\ref{3.20}) and (\ref{3.21}) give a straightforward way
of computing all quantities of weighted potential theory,
necessary for the applications of our results from Section 2.

We now proceed in the same fashion as in the case of two factors,
by defining the weight
$$w(x)=\left(x^{2\alpha_1}|4x-1|^{\alpha_2}|5x-1|^{2\alpha_3}
\right)^{1/(1-2\alpha_1-\alpha_2-2\alpha_3)}, \quad x\in
[0,1/4],$$ where
$$ (\alpha_1,\alpha_2,\alpha_3) \in T:=\{2\alpha_1+\alpha_2+
2\alpha_3 < 1\} \cap \R^3_+.$$ Theorem \ref{thm2.6} applies here
with $\zeta_1=0,\ \zeta_2=1/4$ and $\zeta_3=1/5,$ so that we have
the corresponding values $q_1=1,\ q_2=4$ and $q_3=5$ for
(\ref{2.15}). Each $q_i$ gives a lower bound
\begin{equation} \label{3.22}
l_i(\alpha_1,\alpha_2,\alpha_3):=q_i^{2\alpha_1+\alpha_2+2\alpha_3-1}
e^{(2\alpha_1+\alpha_2+2\alpha_3-1) (F_w-U^{\mu_w}(\zeta_i))},
\ i=1,2,3,
\end{equation}
by (\ref{2.15}), so that
$$t_{\Z}([0,1/4])\ge\inf_{T} \max_i l_i(\alpha_1,\alpha_2,\alpha_3) >
0.1775.$$ This numerical lower bound was found by using
(\ref{3.20}) and a simple C code for generating weighted Leja
points. It should be noted that the most time consuming part of
the computation is finding of the above inf, which is done by a
search over a discrete lattice in $T$. The lower bound of Theorem
\ref{thm3.1} follows at once from (\ref{3.1}).

Using the upper bound $M=0.179335$ in Corollary \ref{cor2.7},
as in the two-factor case, and taking advantage of the ready
numerical results on computing (\ref{3.22}), we find the
following improved bounds for $\alpha_1,\ \alpha_2$ and
$\alpha_3$.

\begin{thm} \label{thm3.4}
The integer Chebyshev polynomials $\{ Q_n \}_{n =1}^{\infty}$
on $[0,1]$ satisfy
\begin{equation} \label{3.23}
Q_n (x)=(x(1 -x))^{[ \alpha_1 n]} (2x -1)^{[ \alpha_2 n]}
(5x^2-5x+1)^{[ \alpha_3 n]} R_n (x), \quad \mbox{as } n
\rightarrow \infty,
\end{equation}
where
\begin{equation} \label{3.24}
0.31 \leq \alpha_1 \leq 0.34,\ 0.11 \leq \alpha_2 \leq 0.14
\mbox{ and } 0.035 \leq \alpha_3 \leq 0.057,
\end{equation}
and $R_n \in {\mathcal P}_n ({\Z})$, $n \in {\N}$.
\end{thm}

The upper bound in Theorem \ref{thm3.1} was obtained by using
all eight factors of (\ref{3.4}) to find an upper bound for
$t_{\Z}([0,1/4])$, with the help of the weight
$$w(x) = \left( \prod_{i=1}^8 |Q_{m_i,i}(z)|^{\beta_i}
\right)^{1/(1-\sum_{i=1}^8 \beta_i)}.$$ If we choose the following
set of values for $\beta_i$'s
$$(0.625,0.11,0.07,0.0032,0.0302,0.0112,0.0048,0.00094),$$
then the upper bound
$$t_{\Z}([0,1/4]) < 0.1791$$
is easily found from Theorem \ref{thm2.3} and (\ref{3.21}) via
another computation, involving weighted Leja points. This implies
by (\ref{3.1}) that
$$t_{\Z}([0,1]) < 0.4232,$$
as claimed in Theorem \ref{thm3.1}.

One can further improve the numerical results on the lower bound
for $t_{\Z}([0,1])$ and the bounds for $\alpha_i$'s, by
considering four and more factors from (\ref{3.4}). The upper
bound of Theorem \ref{thm3.1} can also be improved by optimizing
the choice of $\beta_i$'s. All the details for computations and
suggested improvements of numerical results will be published
separately.

\section{Weighted capacity and potentials}

We give a brief description of the basic facts from the
potential theory with external fields, or weighted potential
theory, for the convenience of the reader. One should consult
Saff and Totik \cite{ST97} for a complete exposition including
the history of this subject.

With ${\mathcal M}(E)$ denoting the class of all positive Borel
measures $\mu$ on ${\C}$ such that $\mu ({\C})=1$ and
supp~$\!\mu \subset E$, consider the following weighted energy
problem (cf. \cite[Section I.1]{ST97}):

\indent For the weighted energy integral
\begin{equation} \label{4.1} I_{w} (\mu):= \int\!\!\int \log
\dis\frac{1}{|z-t|w(z)w(t)} \ d \mu(z)d \mu(t), \quad \mu \in
{\mathcal{M}} (E),
\end{equation}
find
\begin{equation} \label{4.2}
V_{w}:= \dis\inf_{\mu \in {\mathcal M}(E)} I_{w}(\mu),
\end{equation}
and identify the extremal measures, if the infimum in (\ref{4.2})
is attained.

The following is a special case of Theorem I.1.3 in \cite{ST97}.

\begin{prop} \label{prop4.1}
Let $w:E\to[0,+\infty)$ be a continuous function on a compact
set $E \subset \C$ such that $cap(\{z\in E: w(z)>0\})>0$. Then \\
{\rm (a)}   $V_{w}$ of (\ref{4.2}) is finite;\\
{\rm (b)}   there exists a unique $\mu_{w} \in
{\mathcal{M}}(E)$ such that
$I_{w} (\mu_{w}) = V_{w}$;\\
{\rm (c)}   $U^{\mu_{w}}(z)-\log w(z) \geq F_{w}$, for quasi
every
$z \in E$;\\
{\rm (d)}   $U^{\mu_{w}}(z)-\log w(z) \leq F_{w}, \ z \in
S_w$,\\
where $S_w:={\rm supp}\! \ \mu_{w}$ and $F_{w}:=V_{w} + \dis\int
\log w(t) d\mu_{w}(t)$.
\end{prop}

By saying in (c) that a property holds quasi everywhere (q.e.),
we mean that it holds everywhere,  with the possible exception
of a set of zero logarithmic capacity (cf. \cite[Sec.
I.1]{ST97}). The weighted capacity of $E$ is then defined by
\begin{equation} \label{4.3}
cap(E,w):=e^{-V_w}.
\end{equation}
In the case $cap(\{z\in E: w(z)>0\})=0$, we set $cap(E,w)=0.$

It will become clear from the proofs that the $n$-th root
asymptotic behavior of integer Chebyshev polynomials is
essentially equivalent to that of the weighted polynomials, for
$w$ given by (\ref{2.3})-(\ref{2.5}). Therefore, weighted
potential theory provides useful tools for the study of integer
Chebyshev problem, such as the following proposition (see Theorem
III.2.1 and Corollary III.2.6 in \cite{ST97}).

\begin{prop} \label{prop4.2}
Suppose that the assumptions of Proposition \ref{prop4.1} are
satisfied. Then, for any polynomial $P_{n} \in {\mathcal
P}_n(\C)$, we have
\begin{equation} \label{4.4}
|w^{n}(z)P_{n}(z)| \leq \|w^{n}P_{n}\|_{S_{w}} \exp (n (F_{w}-
U^{\mu_{w}}(z) + \log w(z))), \quad z \in \C.
\end{equation}

Assume further that for every point $z \in E$, the set
$\{t:|t-z| < \delta, \ t \in E\}$ has positive capacity for any
$\delta > 0$. Then
\begin{equation} \label{4.5}
\|w^{n}P_{n}\|_{E} = \|w^{n} P_{n}\|_{S_{w}}.
\end{equation}
\end{prop}

\section{Proofs}

The following proof is found in Trigub \cite[p. 316]{Tr71}.
\begin{proof}[Proof of Remark \ref{rem1.1}]
Consider the Chebyshev polynomials for $[-2,2]$, given by
\begin{equation} \label{5.1}
t_n(x)=2\cos(n\arccos(x/2))=2^{-n}\left((x+\sqrt{x^2-4})^n
+(x-\sqrt{x^2-4})^n\right), \quad n \in \N.
\end{equation}
We already observed in Section 1 that $t_n(x)$ is a monic
polynomial with integer coefficients, whose roots are given by
(\ref{1.4}). Schur showed that if $n=p$ is a prime number, then
$t_p(x)/x$ is irreducible over integers (see \cite[p. 228]{Ri90}).
Hence the numbers
$$2\cos\frac{(2k-1)\pi}{2p},\quad k=1,\ldots,p;\ k\neq
\frac{p+1}{2},$$ form a complete set of $p-1$ conjugate
algebraic integers in $[-2,2]$, for any prime $p$. It is clear
that the corresponding roots $\{b_k\}_{k=1}^{p-1}$ of
$F_{p-1}(x)=t_p(x-m-2)/(x-m-2)$, obtained by shifting the above
set by $m+2$, form a complete set of conjugates on $[m,m+4]$.
Let $Q_{p-2}$ be an integer Chebyshev polynomial of degree
$p-2$ for $[1/(m+4),1/m]$. Note that
$$ Q_{p-2}\left(\frac{1}{b_k}\right)=\frac{\tilde
Q_{p-2}(b_k)}{b_k^{p-2}} \neq 0, \quad k=1,\ldots,p-1,$$ where
$\tilde Q_{p-2} \in {\mathcal P}_{p-2}(\Z).$ Indeed, if $\tilde
Q_{p-2}(b_k)=0$ for just one $k$, then this must be true for
all $k=1,\ldots,p-1,\mbox{ i.e., }\tilde Q_{p-2}\equiv 0.$
Since the product $\prod_{k=1}^{p-1} \tilde Q_{p-2}(b_k)$ is a
symmetric form in $b_k$'s with integer coefficients, it may be
written as a polynomial in the elementary symmetric functions of
$b_k$'s, with integer coefficients, by the fundamental theorem
on symmetric forms. Thus the above product must be a nonzero
integer, so that we have
$$\left|\prod_{k=1}^{p-1} Q_{p-2}\left(\frac{1}{b_k}\right)\right| =
\frac{\dis \left| \prod_{k=1}^{p-1} \tilde Q_{p-2}(b_k)
\right|}{\dis \left|\prod_{k=1}^{p-1} b_k^{p-2}\right|} \ge
\frac{1}{\dis \left|\prod_{k=1}^{p-1} b_k\right|^{p-2}}$$ and
\begin{equation} \label{5.2}
\norm{Q_{p-2}}_{[1/(m+4),m]}^{p-1} \ge \frac{1}{\dis
\left|\prod_{k=1}^{p-1} b_k\right|^{p-2}}.
\end{equation}
Observe that
$$ \left|\prod_{k=1}^{p-1} b_k\right| = |F_{p-1}(0)| =
\left|\frac{t_p(m+2)}{m+2}\right| \le
\frac{(m+2+\sqrt{(m+2)^2-4})^p}{2^{p-1}(m+2)}, $$ where the last
inequality follows from (\ref{5.1}). Combining (\ref{5.2}) with
the above estimate and (\ref{1.6}), we obtain that
$$ t_{\Z}([1/(m+4),1/m]) \ge \frac{2}{m+2+\sqrt{(m+2)^2-4}} >
\frac{1}{m+2}.$$
\end{proof}

\begin{proof}[Proof of Proposition \ref{prop1.2}]
For the sequence of polynomials $V_m^k(z),\ k \in \N,$ we have
$$t_{\Z}(L_r) \le \lim_{k \to \infty} \norm{V_m^k}_{L_r}^{km} =
r^{1/m}.$$ Thus the upper bound in (\ref{1.13}) follows. Suppose
that $P_l \in {\mathcal P}_{l}(\Z)$ has a leading coefficient
$b_l \neq 0$. Then we estimate
$$\norm{P_l}_{L_r}^{1/l} = |b_l|^{1/l}\norm{z^l+\ldots}_{L_r}^{1/l}
\ge cap(L_r) = (r/|a_m|)^{1/m},$$ by \cite[p. 155 and p.
135]{Ra95}, which gives the lower bound of (\ref{1.13}).
\end{proof}

\begin{proof}[Proof of Corollary \ref{cor1.3}]
Since $a_m=1$, (\ref{1.14}) follows at once from (\ref{1.13}).
Furthermore, if $P_{l} \in {\mathcal P}_{km}(\Z)$ is of exact
degree $l$, with the leading coefficient $b_l \neq 0$, then
$$\norm{P_{l}}_{L_r} = |b_l|\norm{z^l+\ldots}_{L_r}
\ge (cap(L_r))^l = r^{l/m} \ge r^k = \norm{V_m^k}_{L_r},$$ where
we used \cite[p. 155]{Ra95}. Hence $V_m^k(z)$ is an integer
Chebyshev polynomial of degree $km$ on $L_r,\ k \in \N.$
\end{proof}

\begin{proof}[Proof of Theorem \ref{thm1.4}]
If $V_m^k(z)$ is an integer Chebyshev polynomial of degree $km$ on
$L_r$, for any $k \in \N,$ then (\ref{1.15}) is immediate from the
definition (\ref{1.6}). Therefore, we only need to prove the
second statement of the theorem. It is trivial for $r=0$, so that
we assume $r\in (0,1/|a_m|].$ Suppose to the contrary that there
exists a polynomial $P_{l} \in {\mathcal P}_{km}(\Z)$, of exact
degree $l$, such that
$$\norm{P_{l}}_{L_r} < \norm{V_m^k}_{L_r} = r^k.$$
Let $z_i,\ i=1,\ldots,m,$ be the zeros of $V_m$. Clearly, all
$z_i$'s are inside $L_r$, so that we have by the maximum principle
$$ |P_l(z_i)| \le \norm{P_{l}}_{L_r} < r^k,\quad i=1,\ldots,m.$$
Using a known argument based on the fundamental theorem of
symmetric forms (see Lemma in \cite[p. 181]{Mo94}), we obtain that
$$N=a_m^l \prod_{i=1}^m P_l(z_i) \in \Z.$$ On the other hand, we
estimate
$$\left|a_m^l \prod_{i=1}^m P_l(z_i)\right| \le |a_m|^l
\norm{P_{l}}_{L_r}^m < |a_m|^{km} r^{km} \le 1.$$ Consequently,
this integer N is equal to zero, which means that $P_l(z_i)=0$ for
some $i$. But then the irreducible polynomial $V_m$ must divide
$P_l$.

Assume that $P_l(z)=V_m^d(z) R(z)$, where $d \in \N$ and $R \in
{\mathcal P}_{m(k-d)}(\Z)$, of exact degree $l-md$, does not have
$V_m$ as a factor. It follows that
$$|R(z)|=|P_l(z)|/|V_m(z)|^d<r^{k-d}, \quad z \in L_r,$$
and
$$\norm{R}_{L_r} < r^{k-d}.$$ Hence we can use the same argument for
$R$, to conclude that
$$a_m^{l-md} \prod_{i=1}^m R(z_i) = 0.$$
This implies that $V_m$ divides $R$, contradicting our assumption.
\end{proof}

Before giving the proof of Theorem \ref{thm1.5}, we need to state
two lemmas. The first one shows that if a sequence of polynomials
is composed of only finitely many factors, then the $n$-th root
behavior of this sequence can be essentially described by a fixed
``polynomial-power" function.

\begin{lem} \label{lem5.1}
Suppose that all polynomials $P_n \in {\mathcal P}_n(\C),$ of
exact degrees $n\in \N,$ have finitely many distinct factors
$P_{m_i,i} \in {\mathcal P}_{m_i}(\C),\ i=1,\ldots,K.$ If, for a
compact set $E \subset \C,$
$$ \lim_{n \to \infty} \norm{P_n}_E^{1/n} = A,$$
then there exist $\alpha_i \in (0,1],\ i=1,\ldots,k\le K,$ such
that
$$ \norm{\prod_{i=1}^k |P_{m_i,i}|^{\alpha_i}}_E = A,$$
where $ \sum_{i=1}^k \alpha_i m_i =1.$
\end{lem}

\begin{proof}
We begin by choosing an increasing subsequence $n_j \in \N, \
j=1,2,\ldots,$ such that
$$ \lim_{j\to \infty} \frac{l_i(n_j)}{n_j} =: \alpha_i, \quad
i=1,\ldots,K,$$ where $l_i(n_j)$ is the power of the factor
$P_{m_i,i}$ in $P_{n_j}$. Clearly, $0\le \alpha_i\le 1,\
i=1,\ldots,K,$ and $ \sum_{i=1}^K \alpha_i m_i =1.$ We may assume
that
$$ 0<\alpha_i\le 1,\ i=1,\ldots,k,\quad \mbox{and} \quad \alpha_i=0,\
i=k+1,\ldots,K.$$ Our goal is to show that the factors with
$\alpha_i=0$ do not have influence on the $n$-th root behavior
for the norms of the sequence. If $z$ is not a zero of
$P_{m_i,i},\ i=k+1,\ldots,K,$ then
$$ \lim_{j \to \infty} |P_{n_j}(z)|^{1/n_j} = \lim_{j \to \infty}
\prod_{i=1}^K |P_{m_i,i}(z)|^{l_i(n_j)/n_j} = \prod_{i=1}^k
|P_{m_i,i}(z)|^{\alpha_i}, $$ where convergence in the above
equation is uniform on compact subsets of $\C \setminus \{z:
P_{m_i,i}(z)=0,\ i=k+1,\ldots,K\}.$ Hence
\begin{equation} \label{5.3}
\prod_{i=1}^k |P_{m_i,i}(z)|^{\alpha_i} \le A,
\end{equation}
for any $z \in E$, with finitely many exceptions. But the function
on the left of (\ref{5.3}) is continuous, so that
$$\norm{\prod_{i=1}^k |P_{m_i,i}|^{\alpha_i}}_E \le A.$$
It is easy to obtain the opposite inequality from
\begin{eqnarray*}
A=\lim_{j \to \infty} \norm{P_{n_j}}_E^{1/n_j} &\le& \lim_{j \to
\infty} \norm{\prod_{i=1}^k |P_{m_i,i}|^{l_i(n_j)/n_j}}_E \
\lim_{j \to \infty} \prod_{i=k+1}^K
\norm{P_{m_i,i}}_E^{l_i(n_j)/n_j} \\ &=& \norm{\prod_{i=1}^k
|P_{m_i,i}|^{\alpha_i}}_E.
\end{eqnarray*}
\end{proof}

The following fact is intuitively obvious.
\begin{lem} \label{lem5.2}
Assume that $P_{m_i,i} \in {\mathcal P}_{m_i}(\Z),\ i=1,\ldots,k.$
For any $A>0$ and any set of exponents $\alpha_i > 0,\
i=1,\ldots,k,$ the equation
$$ \prod_{i=1}^k |P_{m_i,i}(x)|^{\alpha_i} = A$$
has only finitely many solutions on the real line.
\end{lem}

\begin{proof}
Suppose that this is not the case, and there are infinitely many
real solutions of the above equation. Since the function $f(x)$ on
the left hand side grows indefinitely, as $x\to \pm \infty$, all
these solutions are contained in a bounded open interval $I$. Note
that in this case $f(x)$ must have infinitely many points of local
maximum in $I,$  with at least one point of accumulation $x_0 \in
I$. Let $x_n \in I,\ n \in \N,$ be the sequence of maxima for
$f(x)$ in $I$, such that $\dis\lim_{n\to \infty} x_n = x_0$ and
$f(x_n)\ge A,\ n\in\N.$ Observing that $f(x_0)\ge A>0$, we
conclude that there exists a 2-dimensional neighborhood $\Delta$
of $x_0$, free of zeros of $P_{m_i,i},\ i=1,\ldots,k.$ Hence
$f(z)$ can be defined as a single valued analytic function in
$\Delta$, which is real valued on $\Delta \cap \R,$ by an
appropriate choice of branches for the powers $\alpha_i.$ It
follows that $f'(x_n)=0,\ n\in\N,$ where $f'(z)$ is also analytic
in $\Delta$. Thus the zeros of $f'(z)$ have a point of
accumulation in its domain of analyticity, forcing this function
to vanish identically in $\Delta$. This implies that $f(z) \equiv
A,\ z \in \Delta$, which can be extended to the whole domain $G$
of definition for $f(z)$. But that gives an immediate
contradiction, as $G$ has the zeros of $f$ on the boundary.
\end{proof}

\begin{proof}[Proof of Theorem \ref{thm1.5}]
We first note that the actual degrees of integer Chebyshev
polynomials on $E$ cannot be bounded, for our assumption that $E$
has infinitely many points would give at once that $cap(E)=1$, by
(\ref{1.5})-(\ref{1.6}). Suppose to the contrary that there are
only finitely many polynomials, with integer coefficients, that
can be factors of $Q_n,\ n\in \N.$ Then we have by Lemma
\ref{lem5.1} that, for a subsequence $\{n_j\}_{j=1}^{\infty}
\subset \N,$
$$ t_{\Z}(E) = \lim_{n_j \to \infty}
\norm{Q_{n_j}}_{E}^{1/n_j} = \norm{\prod_{i=1}^k
|Q_{m_i,i}|^{\alpha_i}}_{E}, $$ where $Q_{m_i,i} \in {\mathcal
P}_{m_i}(\Z),\ 0<\alpha_i\le 1,\ i=1,\ldots,k,$ and $\sum_{i=1}^k
\alpha_i m_i =1.$ Observe that there are only finitely many points
$x_j\in E$, where the function $\prod_{i=1}^k
|Q_{m_i,i}(x)|^{\alpha_i}$ attains its norm on $E$:
$$\prod_{i=1}^k |Q_{m_i,i}(x_j)|^{\alpha_i} = t_{\Z}(E),\
j=1,\ldots,M,$$ according to Lemma \ref{lem5.2}. Let
$U(\delta):=\cup_{j=1}^M\{x\in E:|x-x_j|\le\delta\}.$ By choosing
$\delta>0$ sufficiently small, we can make the logarithmic
capacity of $U(\delta)$ as small as we wish (see Theorem 5.1.4(a)
in \cite[p. 130]{Ra95}). This implies that $t_{\Z}(U(\delta))$ can
also be made arbitrarily small by (\ref{1.10}). In particular, we
can find an integer Chebyshev polynomial $P_l$ for $U(\delta)$
such that
$$ \norm{P_l}_{U(\delta)}^{1/l} < t_{\Z}(E).$$
It follows that
$$\norm{|P_l|^{\eps/l} \prod_{i=1}^k
|Q_{m_i,i}|^{(1-\eps)\alpha_i}}_{U(\delta)} \le
\norm{P_l}^{\eps/l}_{U(\delta)} \left( t_{\Z}(E) \right)^{1-\eps}
<  t_{\Z}(E),$$ for any $\eps \in (0,1).$ Note that
$$ \prod_{i=1}^k |Q_{m_i,i}(x)|^{\alpha_i} < t_{\Z}(E) -
c(\delta),\quad x\in E \setminus U(\delta),$$ where $c(\delta)>0.$
Hence we can estimate for $x\in E \setminus U(\delta)$ that
\begin{eqnarray*}
|P_l(x)|^{\eps/l} \prod_{i=1}^k |Q_{m_i,i}(x)|^{(1-\eps)\alpha_i}
&<& \norm{P_l}^{\eps/l}_{E} \left(t_{\Z}(E) -
c(\delta)\right)^{1-\eps} \\ &=& \left(\frac{
\norm{P_l}^{1/l}_{E}}{t_{\Z}(E) - c(\delta)}\right)^{\eps}
(t_{\Z}(E) - c(\delta)).
\end{eqnarray*}
It is clear that we can now choose $\eps>0$ sufficiently small,
to insure that
$$\norm{|P_l|^{\eps/l} \prod_{i=1}^k
|Q_{m_i,i}|^{(1-\eps)\alpha_i}}_{E} < t_{\Z}(E).$$ But this
immediately implies that the polynomials
$$P_l^{[n\eps/l]} \prod_{i=1}^k Q_{m_i,i}^{[(1-\eps)\alpha_i
n]}, \quad n \in \N,$$ have smaller sup norms on $E$ than those of
integer Chebyshev polynomials, as $n \to \infty$. This is an
obvious contradiction.
\end{proof}
One can generalize Theorem \ref{thm1.5} to certain classes of
compact sets $E\subset\C$. The major element needed in the proof
is that the intersection of $E$ and any lemniscate defined by
$\dis \prod_{i=1}^k |Q_{m_i,i}(z)|^{\alpha_i}=t_{\Z}(E)$ has
integer Chebyshev constant less than $t_{\Z}(E)$.

We are now passing to the proof of Theorem \ref{thm1.6} and
stating an auxiliary result.

\begin{lem} \label{lem5.5}
Let $E \subset \C$ be a compact set. Define $$ E_{\delta}:=\{ z:
|z-w|\le \delta,\ w \in E \}.$$ Then for any $\eps>0$ we can find
$\delta>0$ such that
$$ t_{\Z}(E_{\delta}) - t_{\Z}(E) \le \eps.$$
\end{lem}

\begin{proof}
Take $\eps>0$ and choose $n$ such that $\norm{Q_n}_E^{1/n} \le
t_{\Z}(E)+\eps/2.$ Clearly, $E\subset H:=\{z: |Q_n(z)|_E^{1/n} \le
t_{\Z}(E)+\eps/2\}.$ On the other hand, $H$ is inside the
lemniscate $L_{\eps}:=\{z: |Q_n(z)|_E^{1/n} = t_{\Z}(E)+\eps\},$
by the maximum principle, so that we can set
$\delta:=dist(H,L_{\eps})>0.$ Hence $E_{\delta}$ lies interior to
$L_{\eps},$ and
$$ t_{\Z}(E_{\delta}) \le t_{\Z}(L_{\eps}) \le t_{\Z}(E)+\eps.$$
The last inequality follows by considering a sequence of
polynomials $(Q_n)^m,\ m\in\N,$ on $L_{\eps}.$
\end{proof}
The result of Lemma \ref{lem5.5} can be quantified, provided we
have some knowledge of the geometric properties for $E$. In fact,
one can show that if $E$ consists of finitely many non-degenerate
continua, then
$$t_{\Z}(E_{\delta}) - t_{\Z}(E) \le C(E) \sqrt{\delta},$$
where $C(E)>0$ depends only on $E$.

\begin{proof}[Proof of Theorem \ref{thm1.6}]
We first assume that $cap(Z)=0$, and then obtain a contradiction,
to prove (\ref{1.24}). Let $\tilde{Z}$ be the closure of all zeros
of the integer Chebyshev polynomials $Q_n,\ n\in\N,$ for $E.$
Since the sets $Z$ and $\tilde{Z}$ differ only by countably many
isolated points, we have that
\begin{equation} \label{5.4}
cap(\tilde Z) = cap(Z) = 0
\end{equation}
(cf. Theorem 5.1.4 in \cite[p. 130]{Ra95}). Hence $\Omega:=
\C\setminus{\tilde Z}$ is a connected open set, which follows from
Theorem 5.3.2(a) of \cite[p. 138]{Ra95}. Consider a sequence of
functions
$$ u_n(z):=\frac{1}{n}\log|Q_n(z)|$$
that are subharmonic in $\C$ and harmonic in $\Omega.$ It follows
from Bernstein-Walsh lemma (see, e.g., Theorem 5.5.7(a) in
\cite[p. 156]{Ra95}) that this sequence is bounded on compact
subsets of $\C$. Therefore, we can select a subsequence
$u_{n_j}(z)$, converging to a harmonic function $u(z)$ locally
uniformly in $\Omega.$ Note that $u(z) = u(\bar{z}),\ z\in\Omega,$
which is inherited from the polynomials $Q_{n},\ n\in\N.$ Also,
$$ u(z) \le \log t_{\Z}(E), \quad z \in E\setminus{\tilde Z}.$$
We want to show that all accumulation points for the solutions of
the equation
\begin{equation} \label{5.5}
u(x) = \log t_{\Z}(E), \quad x \in E,
\end{equation}
belong to $E\cap{\tilde Z}.$ Indeed, if $x_0 \in E\setminus{\tilde
Z}$ is such a point, then it must also be a point of accumulation
for the local maxima of $u(x)$ on $\R\cap\Omega.$ Let $x_k,\
k\in\N,$ be those maxima of $u(x)$ such that $\dis \lim_{k\to
\infty} x_k = x_0.$ Consider a 2-dimensional neighborhood
$\Delta\subset\Omega$ of $x_0$. We can define an analytic
completion of $u(z)$ in $\Delta$, denoted by $f(z),$ such that
$\Im f(x_0)=0.$ It is easy to see from the Schwarz integral
formula that $\overline{f(z)} = f(\bar{z}),\ z\in\Delta,$ because
$u(z) = u(\bar{z}),\ z\in\Delta.$ Hence $\Im f(z) = 0,\ z\in
\Delta\cap\R,$ which means that
$$ u'(x_k)=0\ \Rightarrow\ f'(x_k)=0,\quad k \in\N,$$
where by $f'(z)$ we understand the complex derivative of $f(z)$.
It follows that $f'(z)$ vanishes identically in $\Delta$, so that
$f(z)$ and $u(z)$ are identically constant in $\Delta$. But then
$u(z)$ is identically constant in the whole domain $\Omega$, which
cannot be true, because $\Omega$ contains compact sets $H$ of
arbitrarily large capacity and $\norm{u_n}_H \ge \log cap(H)$, for
any $n\in\N$ (cf. Theorem 5.5.4(a) of \cite[p. 155]{Ra95}).

Thus the set $M$ of solutions for (\ref{5.5}) in
$E\setminus{\tilde Z}$ consists of isolated points, i.e., $M$ is
countable and $cap(M)=0$. Furthermore,
$$cap\left((E\cap{\tilde Z})\cup M\right) = 0,$$
by (\ref{5.4}) and Theorem 5.1.4 of \cite[p. 130]{Ra95}. Set
$$ U(\delta):=\{ y\in E: |y-x|\le \delta,\ x \in
(E\cap{\tilde Z})\cup M \}.$$ We choose a sufficiently small
$\delta>0$, so that
$$t_{\Z}(U(\delta)) < t_{\Z}(E),$$
by Lemma \ref{lem5.5} and (\ref{1.10}). Hence there exists
$c_1(\delta)>0$, such that
$$\norm{P_l}_{U(\delta)}^{1/l} < t_{\Z}(E) - c_1(\delta),$$
for integer Chebyshev polynomials $P_l$ on $U(\delta)$ of degree
$l\ge l_0$. Recall that for the integer Chebyshev polynomials
$Q_{n_j}$ on $E$, we have
$$\norm{Q_{n_j}}_E^{1/n_j} < t_{\Z}(E) + \eps_j,\quad j\in\N,$$
where $\lim_{j\to\infty} \eps_j = 0.$ We now let $l_j=l_0+[n_j
\sqrt{\eps_j}]$ and consider sequences of polynomials $\{P_{l_j}^m
Q_{n_j}^m\}_{m=1}^{\infty},\ j\in\N$. Using two preceding
estimates and Young's inequality, we obtain that
\begin{eqnarray} \label{5.6}
\norm{P_{l_j}^m Q_{n_j}^m}_{U(\delta)}^{1/(m(l_j+n_j))} &\le&
\norm{P_{l_j}}_{U(\delta)}^{1/(l_j+n_j)}
\norm{Q_{n_j}}_{U(\delta)}^{1/(l_j+n_j)} \nonumber \\ &<&
(t_{\Z}(E) - c_1(\delta))^{l_j/(l_j+n_j)} (t_{\Z}(E) +
\eps_j)^{n_j/(l_j+n_j)} \nonumber \\ &\le&
\frac{l_j}{l_j+n_j}\left(t_{\Z}(E) - c_1(\delta)\right) +
\frac{n_j}{l_j+n_j}\left(t_{\Z}(E) + \eps_j\right) \nonumber \\
&=& t_{\Z}(E) - \frac{l_j c_1(\delta) - n_j \eps_j}{l_j+n_j} <
t_{\Z}(E),
\end{eqnarray}
for all large $j\in\N.$

Observe that we can find $c_2(\delta)>0$, so that
$$ u(x) < \log\left(t_{\Z}(E) - 2 c_2(\delta)\right),\quad x\in
E\setminus U(\delta),$$ by our construction of the set
$U(\delta).$ Therefore,
$$\norm{Q_{n_j}}_{E\setminus U(\delta)}^{1/n_j} < t_{\Z}(E) -
c_2(\delta),$$ for all sufficiently large $j\in\N.$ This gives the
following estimate
\begin{eqnarray} \label{5.7}
\norm{P_{l_j}^m Q_{n_j}^m}_{E\setminus U(\delta)}^{1/(m(l_j+n_j))}
&\le& \norm{P_{l_j}}_{E}^{1/(l_j+n_j)} (t_{\Z}(E) -
c_2(\delta))^{n_j/(l_j+n_j)} \nonumber \\ &=&
\left(\frac{\norm{P_{l_j}}_{E}^{1/l_j}}{t_{\Z}(E) - c_2(\delta)}
\right)^ \frac{l_j}{l_j+n_j} (t_{\Z}(E)-c_2(\delta)) < t_{\Z}(E),
\end{eqnarray}
where $j$ is selected to be sufficiently large. The last
inequality in (\ref{5.7}) follows because
$\norm{P_{l_j}}_{E}^{1/l_j} < c_3(\delta),\ j\in\N,$ by
Bernstein-Walsh inequality, and because $\lim_{j\to\infty} l_j/n_j
= 0.$ Finally, we combine (\ref{5.6}) and (\ref{5.7}) to obtain
the contradiction:
$$t_{\Z}(E) \le \lim_{m\to\infty} \norm{P_{l_j}^m Q_{n_j}^m}_{E}
^{1/(m(l_j+n_j))} < t_{\Z}(E).$$
\end{proof}

\begin{proof}[Proof of Theorem \ref{thm2.1}]
Observe that if the set $E'=\{z\in E : w(z)>0 \}$ is finite, then
$t_{\Z}(E,w)=0$. Indeed, we can use the regular integer Chebyshev
polynomials $Q_n,\ n\in\N,$ on $E'$, to find that
$$t_{\Z}(E,w) \le \limsup_{n\to\infty} \norm{w^nQ_n}_E^{1/n} \le
\norm{w}_{E'} \lim_{n\to\infty} \norm{Q_n}_{E'}^{1/n} =
\norm{w}_{E'} t_{\Z}(E').$$ But $cap(E')=0$ in this case, so that
$t_{\Z}(E')=0$ by (\ref{1.10}). Thus, (\ref{2.2}) is trivially
true when $E'$ is finite, and we assume that $E'$ has infinitely
many points for the rest of this proof.

We need to find a sequence of polynomials
$$P_n(z)=\sum_{k=0}^n a_k z^k \in {\mathcal P}_n(\Z), \quad
n\in\N,$$ with small weighted norms $\norm{w^nP_n}_E.$ It is
possible to use the Lagrange interpolation in weighted Fekete
points for this purpose. The weighted Fekete points
$\{\zeta_j\}_{j=0}^n \subset E$ are defined as a set of points
maximizing the absolute value of the ``weighted Vandermonde
determinant" (cf. \cite[p. 143]{ST97})
$$V_w(z_0,\ldots,z_n):=\prod_{0\le i<j\le n}
(z_i-z_j)w(z_i)w(z_j)$$ among all $(n+1)$-tuples $\{z_j\}_{j=0}^n
\subset E$. Note that $w(\zeta_j)\neq 0,\ j=0,\ldots,n,$ and we
obtain from the Lagrange interpolation formula that
\begin{equation} \label{5.50}
w^n(z)P_n(z) = \sum_{i=0}^n w^n(\zeta_i)P_n(\zeta_i) \prod_{j\neq
i} \frac{(z-\zeta_j)w(z)}{(\zeta_i-\zeta_j)w(\zeta_i)}.
\end{equation}
Since
$$|V_w(\zeta_0,\ldots,z,\ldots,\zeta_n)| \le
|V_w(\zeta_0,\ldots,\zeta_i,\ldots,\zeta_n)|, \quad z\in E,$$ for
any $i=0,\ldots,n,$ we have that
$$\left|\prod_{j\neq i} \frac{(z-\zeta_j)w(z)}
{(\zeta_i-\zeta_j)w(\zeta_j)}\right| =
\left|\frac{V_w(\zeta_0,\ldots,z,\ldots,\zeta_n)}
{V_w(\zeta_0,\ldots,\zeta_i,\ldots,\zeta_n)}\right| \le 1,\quad
z\in E.$$ It follows at once from (\ref{5.50}) that
\begin{equation} \label{5.51}
\norm{w^nP_n}_E \le \sum_{i=0}^n |w^n(\zeta_i)P_n(\zeta_i)| \le
(n+1)\max_{0\le i\le n} |w^n(\zeta_i)P_n(\zeta_i)|
\end{equation}
(also see Theorem III.1.12 in \cite{ST97}). Observe that
$$ l_i:=w^n(\zeta_i)P_n(\zeta_i)=\sum_{k=0}^n
w^n(\zeta_i)\zeta_i^k a_k, \quad i=0,\ldots,n,$$ are linear forms
in $\{a_k\}_{k=0}^n$, with real coefficients. Applying Minkowski's
theorem (see \cite[p. 73]{Ca97}), we conclude that there exists a
set of integers $\{a_k\}_{k=0}^n$, not all zero, such that
$$|l_i| \le |\det\left(w^n(\zeta_i)\zeta_i^k\right)_{0\le i,k\le
n}|^{1/(n+1)}.$$ But
$$\det\left(w^n(\zeta_i)\zeta_i^k\right)_{0\le i,k\le n} =
V_w(\zeta_0,\ldots,\zeta_n),$$ so that we can find a sequence
$P_n(z)=\sum_{k=0}^n a_k z^k \not\equiv 0,$ satisfying
$$\norm{w^nP_n}_E \le (n+1)|V_w(\zeta_0,\ldots,\zeta_n)|
^{1/(n+1)},\quad n\in\N.$$ Hence
$$\lim_{n\to \infty} \norm{w^nP_n}_E^{1/n} \le \lim_{n\to \infty}
|V_w(\zeta_0,\ldots,\zeta_n)|
^{\frac{1}{n(n+1)}}=\sqrt{cap(E,w)},$$ by Theorem III.1.3 of
\cite[p. 145]{ST97}.
\end{proof}

\begin{proof}[Proof of Theorem \ref{thm2.3}]
Let $P_n \in {\mathcal P}_n(\Z), \ n\in\N,$ be a sequence
polynomials satisfying
$$\lim_{n\to\infty}\norm{w^nP_n}_E^{1/n}=t_{\Z}(E,w),$$
where $w$ is defined in (\ref{2.3}). We construct the following
new sequence of polynomials with integer coefficients:
$$P_n(z)\prod_{i=1}^k Q_{m_i,i}^{l_i(n)}(z), \quad n\in\N,$$
where $l_i(n) \in \N$ are selected so that
$$ \frac{l_i(n)}{n} \searrow \frac{\alpha_i}{1-\alpha}, \quad
\mbox{ as } n\to \infty, \ i=1,\ldots,k.$$ Using (\ref{1.6}), we
obtain that
\begin{eqnarray*}
t_{\Z}(E) &\le& \limsup_{n\to\infty} \norm{P_n\prod_{i=1}^k
Q_{m_i,i}^{l_i(n)}}_E^{1/\left(n+\sum_{i=1}^k m_i l_i(n)\right)} \\
&\le& \limsup_{n\to\infty} \left(\norm{w^nP_n}_E^{1/n}\right)
^{\frac{n}{n+\sum_{i=1}^k m_i l_i(n)}} \\
&\times& \limsup_{n\to\infty}
\left(\prod_{i=1}^k\norm{Q_{m_i,i}}_E^{\frac{l_i(n)}{n}-
\frac{\alpha_i}{1-\alpha}}\right) ^{\frac{n}{n+\sum_{i=1}^k m_i
l_i(n)}} \le \left(t_{\Z}(E,w)\right)^{1-\alpha},
\end{eqnarray*}
because
$$\limsup_{n\to\infty}
\norm{Q_{m_i,i}}_E^{\frac{l_i(n)}{n}- \frac{\alpha_i}{1-\alpha}}
\le 1, \quad i=1,\ldots,k.$$ It now follows from Theorem
\ref{thm2.1} that
$$t_{\Z}(E) \le \left(t_{\Z}(E,w)\right)^{1-\alpha} \le
\left(cap(E,w)\right)^{(1-\alpha)/2}.$$
\end{proof}

\begin{proof}[Proof of Theorem \ref{thm2.4}]
We need to find a solution of the weighted energy problem on $E$,
corresponding to the weight $w$ of (\ref{2.3})-(\ref{2.5}). It
follows from Theorem I.1.3 of \cite{ST97} that there exists a
weighted equilibrium measure $\mu_w$, whose support is a compact
set $S_w \subset E \setminus \cup_{i=1}^k
\{z_{j,i}\}_{j=1}^{m_i}.$ Let $\delta_z$ be a unit point mass at
$z$. Observe that
\begin{equation} \label{5.52}
Q(z)=-\log w(z)=U^{\nu}(z)-\frac{1}{1-\alpha} \sum_{i=1}^k
\alpha_i\log|a_i|,
\end{equation}
where $U^{\nu}$ is the logarithmic potential of the measure
$$\nu:=\frac{1}{1-\alpha} \sum_{i=1}^k \sum_{j=1}^{m_i} \alpha_i
\delta_{z_{j,i}}.$$ It is clear that $\nu$ is a positive Borel
measure of total mass $\nu(\C)=\alpha/(1-\alpha).$ Let $\hat\nu$
be the balayage of $\nu$ from $\Omega$ onto $S_w$ (see, e.g.,
Section II.4 of \cite{ST97}). Then $\hat\nu$ is a positive Borel
measure of the same mass as $\nu$, which is supported on $S_w$.
Furthermore, we can express $\hat\nu$ via harmonic measures
$$\hat\nu=\frac{1}{1-\alpha}\sum_{i=1}^k \sum_{j=1}^{m_i} \alpha_i
\omega(z_{j,i},\cdot,\Omega)$$ (cf. Appendix A.3 of \cite{ST97}).
The potentials of $\nu$ and $\hat\nu$ are related by the equation
\begin{equation} \label{5.53}
U^{\hat\nu}(z)=U^{\nu}(z)+\int_{\Omega} g_{\Omega}(t,\infty)
d\nu(t),
\end{equation}
which holds quasi everywhere on $S_w$ (see Theorem II.4.4 of
\cite{ST97}). Hence the measure
$$\mu:=\frac{1}{1-\alpha}\omega(\infty,\cdot,\Omega) - \hat\nu =
\frac{1}{1-\alpha}\left(\omega(\infty,\cdot,\Omega) - \sum_{i=1}^k
\sum_{j=1}^{m_i} \alpha_i \omega(z_{j,i},\cdot,\Omega) \right)$$
is a probability measure on $S_w$. Using (\ref{5.52}) and
(\ref{5.53}), we obtain for quasi every $z\in S_w$ that
\begin{eqnarray*}
U^{\mu}(z)+Q(z)&=& \frac{1}{\alpha-1}\log{cap(S_w)}-U^{\hat\nu}(z)
+U^{\nu}(z)-\frac{1}{1-\alpha} \sum_{i=1}^k \alpha_i\log|a_i| \\
&=&  \frac{1}{\alpha-1}\left(\log{cap(S_w)} + \sum_{i=1}^k
\alpha_i\log|a_i|\right) - \int_{\Omega} g_{\Omega}(t,\infty)
d\nu(t) \\ &=& \frac{1}{\alpha-1}\left(\log{cap(S_w)} +
\sum_{i=1}^k \alpha_i\log|a_i| + \sum_{i=1}^k \sum_{j=1}^{m_i}
\alpha_i g_{\Omega}(z_{j,i},\infty) \right).
\end{eqnarray*}
Note that $\mu$ has finite logarithmic energy, since it is
composed of harmonic measures. Thus we can apply Theorem I.3.3 to
prove that $\mu$ is the weighted equilibrium measure $\mu_w$ and
that the associated modified Robin constant $F_w$ is given by
(\ref{2.8}). Equation (\ref{2.7}) expresses $cap(E,w)$ through
$w$, $\mu_w$ and $F_w$ (cf. Section I.6 of \cite{ST97}).

We now obtain (\ref{2.10}) from (\ref{2.7}) by the following
simple manipulation. Recall that
$$ g_{\Omega}(z,\infty) = - \log cap(S_w) -
U^{\omega(\infty,\cdot,\Omega)}(z) = \int \log|z-t|
d\omega(\infty,t,\Omega) - \log cap(S_w)$$ (see Chapter 4 of
\cite{Ra95}). Substituting this relation into (\ref{2.8}), we have
that
\begin{eqnarray*}
F_w &=& \frac{1}{\alpha-1}\left(\log{cap(S_w)} +
(1-\alpha)\int\log w(t) d\omega(\infty,t,\Omega)
-\alpha\log{cap(S_w)} \right) \\ &=& - \log{cap(S_w)} - \int\log
w(t) d\omega(\infty,t,\Omega).
\end{eqnarray*}
Hence (\ref{2.10}) is proved too.
\end{proof}

\begin{lem} \label{lem5.11}
Assume that the integer Chebyshev polynomials of $E,\ cap(E)>0,$
satisfy (\ref{2.11}) and (\ref{2.12}), along a subsequence of $n
\to \infty.$ Then
\begin{equation} \label{5.54}
\limsup_{n\to\infty} \norm{w^{m(n)}R_n}_{S_w}^{1/n} \le t_{\Z}(E),
\quad m(n)=\deg R_n,
\end{equation}
where $w$ is given by
(\ref{2.3})-(\ref{2.5}) and $S_w$ is the support of the weighted
equilibrium measure $\mu_w$, corresponding to the weight $w$.
\end{lem}
\begin{proof}
Observe that the actual degree of $R_n$ is $m(n)=n-\sum_{i=1}^k
m_i l_i(n).$ It follows from Theorem I.1.3 of \cite{ST97} that
there exists a weighted equilibrium measure $\mu_w$, whose support
is a compact set $S_w \subset E \setminus \cup_{i=1}^k
\{z_{j,i}\}_{j=1}^{m_i}.$ Hence the factors $Q_{m_i,i},\
i=1,\ldots,k,$ do not vanish on $S_w.$ Estimating
\begin{eqnarray*}
\norm{w^{m(n)}R_n}_{S_w}^{1/n} &\le& \norm{\left(\prod_{i=1}^k
Q_{m_i,i}^{l_i(n)}\right) R_n}_{S_w}^{1/n} \norm{\prod_{i=1}^k
|Q_{m_i,i}|^{\alpha_i m(n)/(1-\alpha) - l_i(n)}}_{S_w}^{1/n} \\
&\le& \norm{Q_n}_E^{1/n} \prod_{i=1}^k \norm{Q_{m_i,i}}
^{\frac{\alpha_i m(n)}{(1-\alpha)n} - \frac{ l_i(n)}{n}}_{S_w},
\end{eqnarray*}
and noting that
$$\lim_{n\to\infty} \norm{Q_{m_i,i}}^{\frac{\alpha_i m(n)}
{(1-\alpha)n} - \frac{ l_i(n)}{n}}_{S_w} = 1,\quad i=1,\ldots,k,$$
by (\ref{2.12}), we obtain (\ref{5.54}).
\end{proof}

\begin{proof}[Proof of Theorem \ref{thm2.5}]
Clearly, the degrees of integer Chebyshev polynomials on $E$ must
be unbounded, because $cap(E)>0,$ so that the assumptions of this
theorem are valid. Since the leading coefficient of $R_n$ is at
least 1 in absolute value and the degree of $R_n$ is
$m(n)=n-\sum_{i=1}^k m_i l_i(n)$, we obtain from Theorem I.3.6 of
\cite{ST97} that
$$\norm{w^{m(n)}R_n}_{S_w} \ge e^{-m(n) F_w}.$$
Thus (\ref{2.13}) follows by taking the $n$-th
root in the above inequality, and using (\ref{5.54}) together with
(\ref{2.12}), as $n\to \infty.$

Finally, (\ref{2.14}) is a direct consequence of (\ref{2.13}) and
(\ref{2.8}), for $E\subset\R.$
\end{proof}

\begin{proof}[Proof of Theorem \ref{thm2.6}]
Recall that the actual degree of $R_n$ is given by
$m(n)=n-\sum_{i=1}^k m_i l_i(n).$ Since $R_n(\zeta)\neq 0$, we
have that
$$|R_n(\zeta)| \ge \frac{N}{q^{m(n)}} \ge \frac{1}{q^{m(n)}},$$
for a positive integer $N$. On the other hand, Theorem III.2.1 of
\cite{ST97} gives the estimate
\begin{equation} \label{5.55}
|R_n(z)| \le \norm{w^{m(n)}R_n}_{S_w} \exp\left( m(n)
(F_w-U^{\mu_w}(z))\right),\quad z \in \C.
\end{equation}
Consequently,
\begin{eqnarray*}
\norm{w^{m(n)}R_n}_{S_w} &\ge& |R_n(\zeta)| \exp\left( -m(n)
(F_w-U^{\mu_w}(\zeta))\right) \\ &\ge& q^{-m(n)} \exp\left( -m(n)
(F_w-U^{\mu_w}(\zeta))\right).
\end{eqnarray*}
Taking the $n$-th root in the above inequality, and using
(\ref{5.54}) together with (\ref{2.12}), we obtain (\ref{2.15}),
as $n\to \infty.$
\end{proof}

\begin{proof}[Proof of Proposition \ref{prop2.8}]
Consider a sequence of integer Chebyshev polynomials $Q_n,\
n\in\N,$ satisfying (\ref{2.11}) and (\ref{2.12}). It is not
difficult to see that we can assume
$$ \frac{l_i(n)}{n} \searrow \alpha_i, \quad
\mbox{ as } n\to \infty, \ i=1,\ldots,k,$$ while preserving the
property
$$\lim_{n\to\infty} \norm{Q_n}_E^{1/n}=t_{\Z}(E).$$
Since
$$t_{\Z}(E) \le t_{\Z}(E\cup H_{\eps}),$$
we only need to show that, for some $\eps>0$,
$$\limsup_{n\to\infty} \norm{Q_n}_{H_{\eps}}^{1/n} \le t_{\Z}(E).$$
Using the same notations as in the proofs of Theorems \ref{thm2.5}
and \ref{thm2.6}, we have that
\begin{eqnarray*}
\limsup_{n\to\infty} \norm{Q_n}_{H_{\eps}}^{1/n} &\le&
\limsup_{n\to\infty} \norm{w^{m(n)}R_n}_{H_{\eps}}^{1/n}
\limsup_{n\to\infty} \prod_{i=1}^k \norm{Q_{m_i,i}}
^{\frac{l_i(n)}{n}-\frac{\alpha_i m(n)}{(1-\alpha)n}}_{H_{\eps}}
\\ &\le& \limsup_{n\to\infty} \norm{w^{m(n)}R_n}_{H_{\eps}}^{1/n}.
\end{eqnarray*}
We next estimate, multiplying (\ref{5.55}) by $w^{m(n)}(z)$,
\begin{equation} \label{5.56}
|w^{m(n)}(z)R_n(z)| \le \norm{w^{m(n)}R_n}_{S_w} \exp\left( m(n)
(F_w-U^{\mu_w}(z)+\log w(z))\right),\quad z \in \C.
\end{equation}
Recall that $S_w$ is a compact set, $S_w \subset E \setminus
\cup_{i=1}^k \{z_{j,i}\}_{j=1}^{m_i}$. Therefore, the potential
$U^{\mu_w}$ is harmonic and bounded on $H_{\eps}$, if $\eps>0$ is
sufficiently small. On the other hand, we can obviously make $\log
w$ smaller than any negative number on $H_{\eps}$, by choosing
$\eps$ small. It follows that
$$F_w-U^{\mu_w}(z)+\log w(z) \le 0, \quad z \in H_{\eps},$$
for some $\eps$, which further implies that
$$\norm{w^{m(n)}R_n}_{H_{\eps}} \le  \norm{w^{m(n)}R_n}_{S_w},$$
by (\ref{5.56}). Using Lemma \ref{lem5.11}, we now obtain that
$$ \limsup_{n\to\infty} \norm{Q_n}_{H_{\eps}}^{1/n} \le
\limsup_{n\to\infty} \norm{w^{m(n)}R_n}_{H_{\eps}}^{1/n} \le
\limsup_{n\to\infty} \norm{w^{m(n)}R_n}_{S_w}^{1/n} \le
t_{\Z}(E).$$
\end{proof}

\begin{proof}[Proof of Lemma \ref{lem3.2}]
Note that the weight $w$ of (\ref{3.5}) is just a special case of
the Jacobi weights of Examples IV.1.17  and IV.5.2 in \cite{ST97}.
Thus our problem on the interval $[0, 1/4]$ is easily reduced to
that on the interval $[-1,1]$ considered in \cite{ST97}, with the
help of the change of variable $x \rightarrow (x+1)/8$.  We obtain
from Example IV.1.17 of \cite{ST97} that $S_w = [a,b]$, with $a$
and $b$ given by (\ref{3.6}) (see (1.27) and (1.28) in \cite[p.
207]{ST97}).  Similarly,  Example IV.5.2 gives (\ref{3.8}) after
this change of variable.

Consider the following natural extension for $w(x)$ of
(\ref{3.5}):
$$ w(z) := |z|^{\frac{2\alpha_1}{1-2\alpha_1-\alpha_2}}
|1-4z|^{\frac{\alpha_2}{1-2\alpha_1-\alpha_2}}, \quad z \in
{\C}.$$ It follows from Theorem I.1.3 of \cite{ST97} that
\begin{equation} \label{5.57}
F_w - U^{\mu_w} (z) = - \log w(z),
\end{equation}
for quasi every $z \in [a,b]$ (i.e., with the exception of a set
of zero capacity). Denote the right hand side of (\ref{3.9}) by
$h(z).$ Then $F_w - U^{\mu_w} (z) -h(z)$ is a harmonic function in
$\Omega=\C\setminus [a,b]$, such that
$$F_w - U^{\mu_w}(z)-h(z) =0$$
for quasi every $z \in [a,b] = \partial \Omega$, by (\ref{5.57})
and the basic properties of Green functions (see \cite[p.
14]{Ts75}). Using the uniqueness theorem for the solution of the
Dirichlet problem in $\Omega$ (cf. Theorem III.28 and its
Corollary in \cite{Ts75}), we conclude that
$$F_w - U^{\mu_w} (z) \equiv h(z), \quad z \in {\C}.$$
Thus $F_w$ can be found from
$$F_w=\lim_{z\to\infty} \left( U^{\mu_w}(z) + h(z) \right),$$
or from (\ref{2.8}). We obtain the explicit representation of
(\ref{3.7}) by expressing the Green functions of (\ref{3.9}) via
the conformal mappings of $\Omega$ onto the exterior of the unit
disk. Indeed, introducing these conformal mappings by
$$ \Phi_{\infty} (z) := \frac{2z -a -b +2 \sqrt{ (z -a)(z -b)}}{b -a},
\quad z \in \Omega,$$
$$\Phi_0 (z) := \frac{2z^{-1}-b^{-1}-a^{-1} + 2 \sqrt{(z^{-1}
-b^{-1})(z^{-1} -a^{-1})}}{b^{-1} - a^{-1}}, \quad z\in\Omega,$$
and
\begin{eqnarray*}
\Phi_{1/4} (z) := \frac{2 (z - 1/4)^{-1}-(b-1/4)^{-1} -(a -
1/4)^{-1}} {(b - 1/4)^{-1} - (a - 1/4)^{-1}} + \nonumber \\
\frac{2 \sqrt{((z - 1/4)^{-1} - (a - 1/4)^{-1})((z - 1/4)^{-1}
- (b - 1/4)^{-1})}}{(b - 1/4)^{-1} - (a - 1/4)^{-1}}, \quad z \in \Omega,
\end{eqnarray*}
we observe that
$$\Phi_{\infty} ( \infty ) = \infty, \quad \Phi_0 (0) = \infty \quad
\mbox{ and } \quad \Phi_{1/4} (1/4) = \infty.$$ Hence
\begin{eqnarray*}
g _{\Omega}(z, \infty) = \log | \Phi_{\infty} (z)|, \quad
g_{\Omega}(z,0) = \log | \Phi_0 (z)| \\
\mbox{ and } \quad g_{\Omega} (z, 1/4) = \log |\Phi_{1/4} (z)|,
\quad z \in \Omega,
\end{eqnarray*}
by Theorem I.17 of \cite[p. 18]{Ts75}.
\end{proof}

\begin{proof}[Proof of Theorem \ref{thm3.3}]
This theorem is an immediate application of Corollary \ref{cor2.7}
to the integer Chebyshev polynomials on $[0,1/4]$. Using
(\ref{1.23}) and (\ref{3.1}), we obtain the upper bound
$$t_{\Z}([0,1/4]) \le 0.42347945^2 < 0.179335 = M.$$
Then we choose $\zeta_1=0$ and $\zeta_2=1/4$ to produce the
inequalities, defining the region $G$ of Figure \ref{fig1}, by
(\ref{2.16}). The values of $F_w - U^{\mu_w}(\zeta_i)$ are readily
found from (\ref{3.9}) and the explicit formulas for the Green
functions, obtained in the proof of Lemma \ref{lem3.2}. Figure
\ref{fig1}, as well as the bounds for $\alpha_1$ and $\alpha_2$,
is generated by Matlab.
\end{proof}


\end{document}